\pdfoutput=1
\newcommand{\ExecuteBeforeCleveref}{%
    \usepackage[algo2e,ruled,vlined,linesnumbered]{algorithm2e}%
}

\documentclass{shinyart}

\usepackage[utf8]{inputenc}

\usepackage{shinybib}

\usepackage{autonum}
\usepackage{booktabs}
\usepackage[binary-units]{siunitx}

\title{Efficient high-resolution RF pulse design applied to simultaneous multi-slice excitation}
\author{Christoph Stefan Aigner\thanks{%
    Institute of Medical Engineering, Graz University of Technology, Kronesgasse 5, and BioTechMed Graz, 8010 Graz, Austria (\email{christoph.aigner@tugraz.at})}
    \and Christian Clason\thanks{%
    Faculty of Mathematics, University of Duisburg-Essen, 45117 Essen, Germany}
    \and Armin Rund\thanks{%
    Institute for Mathematics and Scientific Computing, University of Graz, Heinrichstrasse 36, 8010 Graz, Austria}
    \and Rudolf Stollberger\footnotemark[1]
}
\date{November 23, 2015}

\hypersetup{
    pdftitle={Efficient high-resolution RF pulse design applied to simultaneous multi-slice excitation}
    pdfauthor={Christoph Stefan Aigner, Christian Clason, Armin Rund, Rudolf Stollberger},
    pdfkeywords={pulse design, optimal control, second-order methods, simultaneous multi-slice excitation}
} 

\begin{document}
\maketitle

\begin{abstract} 
    RF pulse design via optimal control is typically based on gradient and quasi-Newton approaches and therefore suffers from slow convergence. We present a flexible and highly efficient method that uses exact second-order information within a globally convergent trust-region CG-Newton method to yield an improved convergence rate. The approach is applied to the design of RF pulses for single- and simultaneous multi-slice (SMS) excitation and validated using phantom and in-vivo experiments on a \SI{3}{\tesla} scanner using a modified gradient echo sequence. 

    \vspace*{0.5\baselineskip}
    \noindent{\color{structure}Keywords}\quad pulse design, optimal control, second-order methods, simultaneous multi-slice excitation
\end{abstract}

\section{Introduction}
\label{sec:intro}

For many applications in MRI there is still demand for the optimization of selective RF excitation, e.g., for simultaneous multi-slice excitation \cite{Larkman01, Ugurbil13}, UTE imaging \cite{Bydder12}, or velocity selective excitation \cite{Rochefort06}.
To achieve a well-defined slice profile at high field strength while meeting B$_1$ peak amplitude limitations is a challenge for RF pulse design and becomes especially critical for quantitative methods.
Correspondingly, many approaches for general pulse design have been proposed in the literature.
RF pulses with low flip angles can be designed using the small tip angle simplification \cite{Pauly89}, which makes use of an approximation of the Bloch equation to compute a pulse via the Fourier transform of the desired slice profile. 
However, this simplification breaks down for large flip angles. 
The resulting excitation error for large flip angle pulses can be reduced by applying the Shinnar--Le Roux (SLR) technique \cite{Pauly91} or optimization methods, e.g., simulated annealing, evolutionary approaches or optimal control \cite{Shen01, Wu91, Khalifa01, Conolly86, Xu08, Lapert12, Vinding12}.
The SLR method is based on the hard pulse approximation and a transformation of the excitation problem, allowing to solve the excitation problem recursively by applying fast filter design algorithms such as the Parks--McClellan algorithm \cite{Pauly91}. 
Originally, this approach only covered special pulses such as \ang{90} and \ang{180} excitation or refocusing, but Lee \cite{Lee07} generalized this approach to arbitrary flip angles with an exact parameter relation. 
Despite its limitations due to neglected relaxation terms and sensitivity to B$_1$ inhomogeneities, it found widespread use (see, e.g., \cite{Lee07, Balchandani10, Grissom12, Ma14}) and is considered to be the gold standard for large tip angle pulse design. 
An alternative approach is based on optimizing a suitable functional; see, e.g., 
\cite{Ulloa04, Xu08, Pauly91, Gezelter90, Buonocore93}. In particular, optimal control (OC) approaches involve the solution of the Bloch equation describing the evolution of the magnetization vector in an exterior magnetic field \cite{Conolly86, Khaneja05, Skinner12, Bonnard14, Lapert12, Grissom09, Xu08,Vinding12}. They often lead to better excitation profiles due to a more accurate design model and are increasingly used in MRI, for instance, to perform multidimensional and multichannel RF design \cite{Xu08, Grissom09}, robust 2D spatial selective pulses \cite{Vinding12} and saturation contrast \cite{Lapert12}.
In addition, arbitrary flip angles and target slice profiles, as well as inclusion of additional physical effects such as, e.g., relaxation can be handled.
However, so far OC approaches are limited by the computational effort and require a proper modeling of the objective. In particular, standard gradient-based approaches suffer from slow convergence, imposing significant limitations on the accuracy of the obtained slice profiles.
On the other hand, Newton methods show a locally quadratic convergence, but require second-order information which in general is expensive to compute 
\cite{Anand12}. Approximating the Hessian using finite differences causes loss of quadratic convergence due to the lack of exact second-order information and typically requires significantly more iterations. Superlinear convergence can be obtained using quasi-Newton methods based on exact gradients \cite{Fouquieres11}, although their performance can be sensitive to implementation details.
The purpose of this work is to demonstrate that for the OC approach to pulse design, it is in fact possible to use exact second-order information while avoiding the need of computing the full Hessian, yielding a highly efficient numerical method for the optimal control of the full time-dependent Bloch equation. 
In contrast to \cite{Anand12} (which uses black-box optimization method and symbolically calculated Hessians based on an effective-matrix approximation of the Bloch equation),
we propose a matrix-free Newton--Krylov method \cite{Knoll2004} using first- and second-order derivatives based on the adjoint calculus \cite{Hinze:2001} together with a trust-region globalization \cite{Steihaug83}; for details we refer to \cref{sec:optimization}.
Recently, similar matrix-free Newton--Krylov approaches with line search globalization were presented for optimal control of quantum systems in the context of NMR pulse sequence design \cite{Ciaramella15_1, Ciaramella15_2}.
In comparison, the proposed trust-region framework significantly reduces the computational effort, particularly for the initial steps far away from the optimum.
The effectiveness of the proposed method is demonstrated for the design of pulses for single and simultaneous multi-slice excitation~(SMS).

SMS excitation is increasingly used to accelerate imaging experiments \cite{Larkman01,Nunes06,Feinberg13, Ugurbil13}. 
Conventional design approaches, based on a superposition of phase-shifted sub-pulses \cite{Mueller88} or sinusoidal modulation \cite{Larkman01}, typically result in a linear scaling of the B$_1$ peak amplitude, a quadratic peak power and a linear increase in the overall RF power \cite{Norris11, Auerbach13}. The required maximal B$_1$ peak amplitude of conventional multi-slice pulses therefore easily exceeds the transmit voltage of the RF amplifier. In this case, clipping will occur, while rescaling will decrease and limit the maximal flip angle of such a pulse. On the other hand, restrictions of the specific absorption rate limit the total (integrated) B$_1$ power and therefore the maximal number of slices as well as the pulse duration and flip angle. 
The increase of B$_1$ power can be addressed by the Power Independent of Number of Slices (PINS) technique \cite{Norris11}, which was extended to the kT-PINS method \cite{Sharma13} to account for B$_1$ inhomogeneities. This approach leads to a nearly slice-independent power requirement, but the periodicity of the resulting excitation restricts the slice orientation and positioning.
Furthermore, the slice profile accuracy is reduced \cite{Feinberg13}, and a limited ratio between slice thickness and slice distance may further restrict possible applications.
The combination of PINS with regular multi-band pulses was shown to reduce the overall RF power by up to \SI{50}{\percent} (MultiPINS \cite{Eichner14}) and was applied to refocusing pulses in a multi-band RARE sequence with 13 slices \cite{Gagoski15}.

A different way to reduce the maximum B$_1$ amplitude is to increase the pulse length; however, this stretching increases the minimal echo and repetition times and decreases the RF bandwidth, thus reducing the slice profile accuracy \cite{Auerbach13}.
Applying variable rate selective excitation \cite{Conolly88, Setsompop12-2} avoids this problem but leads to an increased sensitivity to slice profile degradations at off-resonance frequencies. 
In addition, they require specific sequence alterations, e.g., variable slice gradients or gradient blips. 
Instead of using the same phase for all sub-slices, the peak power can be reduced by changing the uniform phase schedule to a different phase for each individual slice \cite{Wong12}.
Alternative approaches\cite{Zhu14, Sharma15} using phase-matched excitation and refocusing pairs show that a nonlinear phase pattern can be corrected by a subsequent refocusing pulse.
Another way to reduce the power deposition and SAR of SMS pulses is to combine them with parallel transmission \cite{Katscher08}. This allows to capitalize transmit sensitivities in the pulse design and leads to a more uniform excitation with an increased power efficiency \cite{Wu13, Poser14}. Recently, Guerin et al. \cite{Guerin14} demonstrated that it is possible to explicitly control both global and local SAR as well as the peak power using a spokes-SMS-pTx pulse design. 

The focus of this work, however, is on single channel imaging, where we apply our OC-based pulse design for efficient SMS pulse optimization using a direct description of the desired magnetization pattern. Its flexible formulation allows a trade-off between the slice profile accuracy and the required pulse power and is well suited for the reduction of power and amplitude requirements of such pulses, even for a large number of slices or large flip angles or in presence of relaxation. The efficient implementation of the proposed method allows to optimize for SMS pulses with a high spatial resolution to achieve accurate excitation profiles. The RF pulses are designed to achieve a uniform effective echo time and phase for each slice and use a constant slice-selective gradient, allowing to insert the RF pulse into existing sequences and opening up a wide range of applications.

\section{Theory} 
\label{sec:theory}

This section is concerned with the description of the optimal control approach to RF pulse design as well as of the proposed numerical solution approach.

\subsection{Optimal control framework}

Our OC approach is based on the full time-dependent Bloch equation, which describes the temporal evolution of the ensemble magnetization vector ${M}(t)=({M_x}(t),{M_y}(t),{M_z}(t))^T$ due to a transient external magnetic field ${B}(t)$ as the solution of the ordinary differential equation (ODE)
\begin{equation}
    \label{eq:bloch1}
    \left\{\begin{aligned}
            \dot{{M}}(t) &= \gamma {B}(t) \times {M}(t) + {R}({M}(t)),\qquad t>0,\\
            {M}(0) &={M^0},  
    \end{aligned}\right.
\end{equation}
where $\gamma$ is the gyromagnetic ratio, ${M^0}$ 
is the initial magnetization and
\begin{equation}
    {R(M}(t)) = (-{M_x}(t)/T_2,-{M_y}(t)/T_2,-({M_z}(t)-M_{0})/T_1)^T
\end{equation}
denotes the relaxation term with relaxation times $T_1,T_2$ and the equilibrium magnetization $M_{0}$. To encode spatial information in MR imaging, the external magnetic field ${B}$ (and thus the magnetization vector) depends on the slice direction $z$, hence the Bloch equation can be considered as a parametrized family of three-dimensional ODEs. In the on-resonance case and ignoring spatial field inhomogeneities, the Bloch equation can be expressed in the rotating frame as 
\begin{equation}
    \label{eq:bloch_rot}
    \left\{\begin{aligned}
            {\dot{M}}(t;z) &= {A}({u}(t);z) {M}(t;z)+ {b}(z), \qquad t>0,\\
            {M}(0;z) &= {M^0}(z),  
    \end{aligned}\right.
\end{equation}
where the control ${u}(t)=({u_x}(t),{u_y}(t))$ describes the RF pulse, 
\begin{equation}
    \label{eq:matrix}
    {A}({u};z) = \begin{pmatrix}
        -\frac1{T_2}       & \gamma {G_z}(t)z    & \gamma {u_y}(t) B_1 \\
        - \gamma {G_z}(t)z  & -\frac1{T_2}       & \gamma {u_x}(t) B_1 \\
        -\gamma {u_y}(t) B_1 & -\gamma {u_x}(t) B_1 & -\frac1{T_1}
    \end{pmatrix}, \qquad
{b}(z) = \begin{pmatrix} 0\\0\\ \frac{{M_{0}}}{T_1} \end{pmatrix},
\end{equation}
and $G_z$ is the slice-selective gradient; see, e.\,g., \cite[Chapter 6.1]{Nishimura96}. 

The OC approach consists in computing for given initial magnetization ${M^0}(z)$ the RF pulse ${u}(t)$, $t\in[0,T_u]$, that minimizes the squared error at read-out time $T>T_u$ between the corresponding solution ${M}(T;z)$ of \eqref{eq:bloch_rot} and a prescribed slice profile ${M_d}(z)$ for all $z\in [-a,a]$ together with a quadratic cost term modeling the SAR of the pulse, i.e., solving
\begin{equation}\label{eq:costfunc}
    \min\limits_{({u},{M}) \text{ satisfying \eqref{eq:bloch_rot}}} \quad J({M},{u}) = \frac{1}{2} \int_{-a}^a |{M}(T;z)-{M_d}(z)|_2^2 \,dz + \frac\alpha2\int_0^{T_u} |{u}(t)|_2^2\,dt.
\end{equation}
The parameter $\alpha>0$ controls the trade-off between the competing goals of slice profile attainment and SAR reduction.

\subsection{Adjoint approach}\label{sec:optimization}

The standard gradient method for solving \eqref{eq:costfunc} consists of computing for given ${u}^k$ the gradient ${g}({u}^k)$ of $j({u}):=J({M}({u}),{u})$ and setting ${u}^{k+1} = {u}^k-s^k {g}({u}^k)$ for some suitable step length $s^k$. The gradient can be calculated efficiently using the adjoint method, which in this case yields
\begin{equation}
    \label{eq:gradient}
    \begin{aligned}[t]
        {g}({u}^k)(t) &= \alpha {u}(t) + \gamma B_1 \begin{pmatrix} 
            \int_{-a}^a {M_z}(t;z){P_y}(t;z) - {M_y}(t;z){P_z}(t;z)\,dz\\[0.5ex]
            \int_{-a}^a {M_z}(t;z){P_x}(t;z) - {M_x}(t;z){P_z}(t;z)\,dz
        \end{pmatrix}\\
        &=: \alpha {u}(t) + 
        \begin{pmatrix}
            \int_{-a}^a {M}(t;z){A_1} {P}(t;z)\,dz\\[0.5ex]
            \int_{-a}^a {M}(t;z){A_2} {P}(t;z)\,dz
        \end{pmatrix},
        \qquad 0\leq t \leq T_u,
    \end{aligned}
\end{equation}
where ${M}$ is the solution to \eqref{eq:bloch_rot} for ${u}={u}^k$ and $0<t\leq T$, ${P}$ is the solution to the adjoint (backward in time) equation
\begin{equation}
    \label{eq:adjoint}
    \left\{\begin{aligned}
            -\dot{{P}}(t;z)&= {A}({u}(t);z)^T {P}(t;z),\qquad 0\leq t< T,\\
            {P}(T;z)&={M}(T;z)-{M_d}(z),
    \end{aligned}\right.
\end{equation}
and for the sake of brevity, we have set
\begin{equation}
{A_1} := \gamma B_1 \begin{pmatrix} 0 & 0 & 0 \\ 0 & 0 & -1 \\ 0 & 1 &0\end{pmatrix},\qquad
{A_2} := \gamma B_1 \begin{pmatrix} 0 & 0 & -1 \\ 0 & 0 & 0 \\ 1 & 0 &0\end{pmatrix}.
\end{equation}

However, this method requires a line search to converge and usually suffers from slow convergence close to a minimizer. 
This is not the case for Newton's method (which is a second-order method and converges locally quadratically), where one additionally computes the Hessian ${H}({u}^k)$ of $j$ at $u^k$, solves for $\delta {u}$ in 
\begin{equation}\label{eq:Newton_step}
    {H}({u}^k)\delta {u} = -{g}({u}^k),
\end{equation}
and sets ${u}^{k+1} = {u}^k+\delta {u}$. While the full Hessian ${H}({u}^k)$ is very expensive to compute in practice, solving \eqref{eq:Newton_step} using a Krylov method such as conjugate gradients (CG) only requires computing the Hessian \emph{action} ${H}({u}^k){h}$ for a given direction ${h}$ per iteration; see, e.g., \cite{Knoll2004}. The crucial observation in our approach is that the adjoint method allows computing this action exactly (e.g., without employing finite difference approximations) and without knowledge of the full Hessian. Since Krylov methods usually converge within very few iterations, this so-called ``matrix-free'' approach amounts to significant computational savings.
To derive a procedure for computing the Hessian action ${H}({u}^k){h}$ for a given direction ${h}$ directly, we start by differentiating \eqref{eq:gradient} with respect to ${u}$ in direction $h$ and applying the product rule. This yields
\begin{equation}
    \label{eq:Hessian}
    [{H}({u}^k){h}](t)  = \alpha {h}(t) +    \begin{pmatrix}
        \int_{-a}^a \delta {M}(t;z){A_1} {P}(t;z) +  {M}(t;z){A_1}\delta {P}(t;z)\,dz\\[0.5ex]
        \int_{-a}^a \delta {M}(t;z){A_2} {P}(t;z) +  {M}(t;z){A_2}\delta {P}(t;z)\,dz
    \end{pmatrix},
    \quad 0\leq t \leq T_u,
\end{equation}
where $\delta {M}$ -- corresponding to the directional derivative of ${M}$ with respect to ${u}$ -- is given by the solution of the linearized state equation
\begin{equation}
    \label{eq:bloch_lin}
    \left\{\begin{aligned}
            \dot{{\delta {M}}}(t;z) &= {A}({u}^k;z) {\delta {M}}(t;z)+ {A}'({h}){M}, \qquad 0<t\leq T,\\
            \delta {M}(0;z) &= (0,0,0)^T,
    \end{aligned}\right.
\end{equation}
with 
\begin{equation}
    {A}'({h})= \gamma B_1
    \begin{pmatrix}
        0       & 0        & {h_y}(t)\\
        0       & 0        & {h_x}(t)\\
        -{h_y}(t) & - {h_x}(t) & 0
    \end{pmatrix},
\end{equation}
and $\delta {P}$ -- corresponding to the directional derivative of ${P}$ with respect to ${u}$ -- is the solution of the linearized adjoint equation
\begin{equation}
    \label{eq:adjoint_lin}
    \left\{\begin{aligned}
            -\dot{\delta {P}}(t;z) &= {A}({u}^k;z)^T {\delta {P}}(t;z)+ {A}'({h})^T {P},\qquad 0\leq t<T,\\
            \delta {P}(T;z) &= \delta {M}(T;z).
    \end{aligned}\right.
\end{equation}
This characterization can be derived using formal Lagrangian calculus and rigorously justified using the implicit function theorem; see, e.g., \cite[Chapter 1.6]{Hinze:2009}.
Since \eqref{eq:Hessian} can be computed by solving the two ODEs \eqref{eq:bloch_lin} and \eqref{eq:adjoint_lin}, the cost of computing a single Hessian action is comparable to that of a gradient evaluation; cf.~\eqref{eq:gradient}. This has already been observed in the context of seismic imaging \cite{Santosa:1988}, meteorology \cite{Wang:1992}, and optimal control of partial differential equations \cite{Hinze:2001}, but has received little attention so far in the context of optimal control of ODEs.

One difficulty is that the Bloch equation \eqref{eq:bloch_rot} is bilinear since it involves the product of the unknowns ${u}$ and ${M}$. Hence, the optimal control problem \eqref{eq:costfunc} is not convex and the Hessian ${H}({u})$ is not necessarily positive definite (or even invertible), thus precluding a direct application of the CG-Newton method. We therefore embed Newton's method into the trust-region framework of Steihaug \cite{Steihaug83}, where a breakdown of the CG method is handled by a trust-region step and the trust region radius is continually adapted. This allows global convergence (i.e., for any starting point) to a local minimizer as well as transition to fast quadratic convergence; see \cite{Steihaug83}.
As an added advantage, computational time is saved since the CG method is usually not fully resolved far away from the optimum. 
The full algorithm is given in \cref{sec:algo}. 

\subsection{Discretization}

For the numerical computation of optimal controls, both the Bloch equation \eqref{eq:bloch_rot} and the optimal control problem in \eqref{eq:costfunc} need to be discretized. 
Here, the time interval $[0,T]$ is replaced by a time grid $0=t_0<\dots<t_N=T$ with time steps $\Delta t_m:=t_m -t_{m-1}$, chosen such that $t_{N_u} = T_u<T$ for some $N_u<N$.
The domain $[-a,a]$ is replaced by a spatial grid $-a = z_1<\dots<z_Z = a$ with grid sizes $\Delta z_m:=z_m-z_{m-1}$. 
We note that for each $z_i$, the corresponding ODEs can be solved independently and in parallel.
The Bloch equation is discretized using a Crank--Nicolson method, where the state ${M}$ is discretized as continuous piecewise linear functions with values ${M}_m:={M}(t_m)$, and the controls ${u}$ are treated as piecewise constant functions, i.e., ${u}=\sum_{m=1}^{N_u} {u_m} \chi_{(t_{m-1},t_m]}(t)$, where $\chi_{(a,b]}$ is the characteristic function of the half-open interval~$(a,b]$.

For the efficient computation of optimal controls, it is crucial that both the gradient and the Hessian action are computed in a manner consistent with the chosen discretization. 
This implies that the adjoint state ${P}$ has to be discretized as piecewise constant using an appropriate time-stepping scheme \cite{BecMeiVex_07}, and that the linearized state $\delta{M}$ and the linearized adjoint state $\delta{P}$ have to be discretized in the same way as the state and adjoint state, respectively.
Furthermore, the conjugate gradient method has to be implemented using the scaled inner product $\langle {u},{v}\rangle := \sum_{m=1}^{N_u}\Delta t_m{u}_m{v}_m$ and the corresponding induced norm $\|{u}\|^2:=\langle {u},{u}\rangle$.
For completeness, the resulting schemes and discrete derivatives are given in \cref{sec:discrete}.

\section{Methods}
\label{sec:methods}

This section describes the computational implementation of the proposed pulse design and the experimental protocol for its validation.

\subsection{Pulse design}
The OC approach described in \cref{sec:theory} is implemented in MATLAB (The MathWorks, Inc., Natick, USA) using the Parallel Toolbox for parallel solution of the (linearized) Bloch and adjoint equation for different values of $z_i$. In the spirit of reproducible research, the code used to generate the results in this paper can be downloaded from \url{https://github.com/chaigner/rfcontrol/releases/v1.2}.

The initial magnetization vector is set to equilibrium, i.e., ${M^0}(z) =M_0(0,0,1)^T$.
The slice-selective gradient $G_z(t)$ is extracted out of a standard Cartesian GRE sequence simulation and consists of a trapezoidal shape of length \SI{2.56}{\milli\second} that is followed by a re-phasing part of length \SI{0.92}{\milli\second} to correct the phase dispersion using the maximal slew rate; i.e., $T_u =$ \SI{2.56}{\milli\second} and $T=$ \SI{3.48}{\milli\second} with a temporal resolution of $\Delta t=$ \SI{5}{\micro\second} for the single-slice excitation (see dashed line in \cref{fig:single:1}) and $T_u =$ \SI{10.24}{\milli\second} and $T=$ \SI{13.92}{\milli\second} with a temporal resolution of $\Delta t=$ \SI{20}{\micro\second} for the SMS excitation (see dashed line in \cref{fig:sms:a2}). This corresponds in both cases to $N=697$ uniform time steps for the time interval $[0,T]$ and $N_u=512$ time steps for the control interval $[0,T_u]$. 
For the spatial computational domain, $a=$ \SI{0.5}{\metre} is chosen to consider typical scanner dimensions; the domain $[-a,a]$ is discretized using $Z=5001$ equidistant points to achieve a homogeneous spatial resolution of $\Delta z=$ \SI{0.2}{\milli\metre}. 

For the desired magnetization vector, we consider three examples:
\begin{figure}
    \centering
    \subcaptionbox{${M_d}$ for single slice (zoom)\label{fig:initial:3_zoom}}{%
        \includegraphics[width=0.3\textwidth]{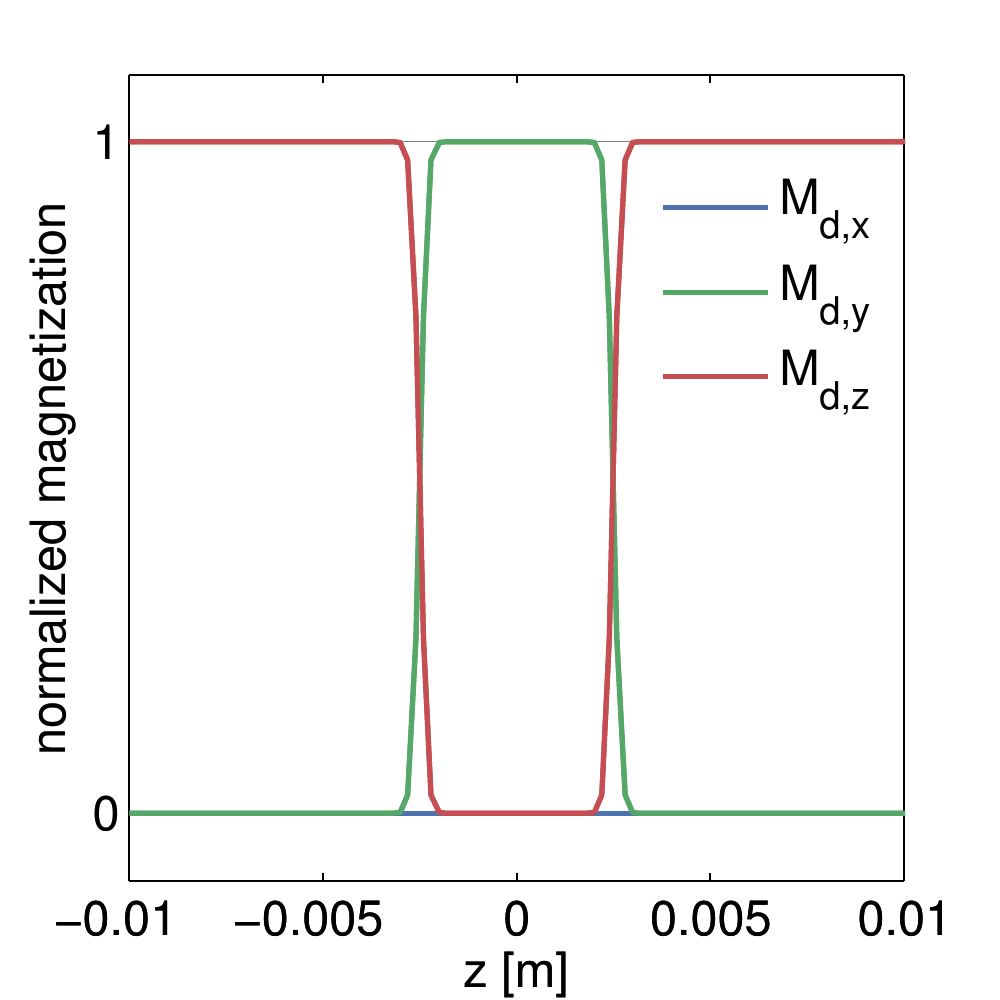}		
    }
    \hfill
    \subcaptionbox{${M_d}$ for SMS 6 (zoom)\label{fig:initial:4_zoom}}{%
        \includegraphics[width=0.3\textwidth]{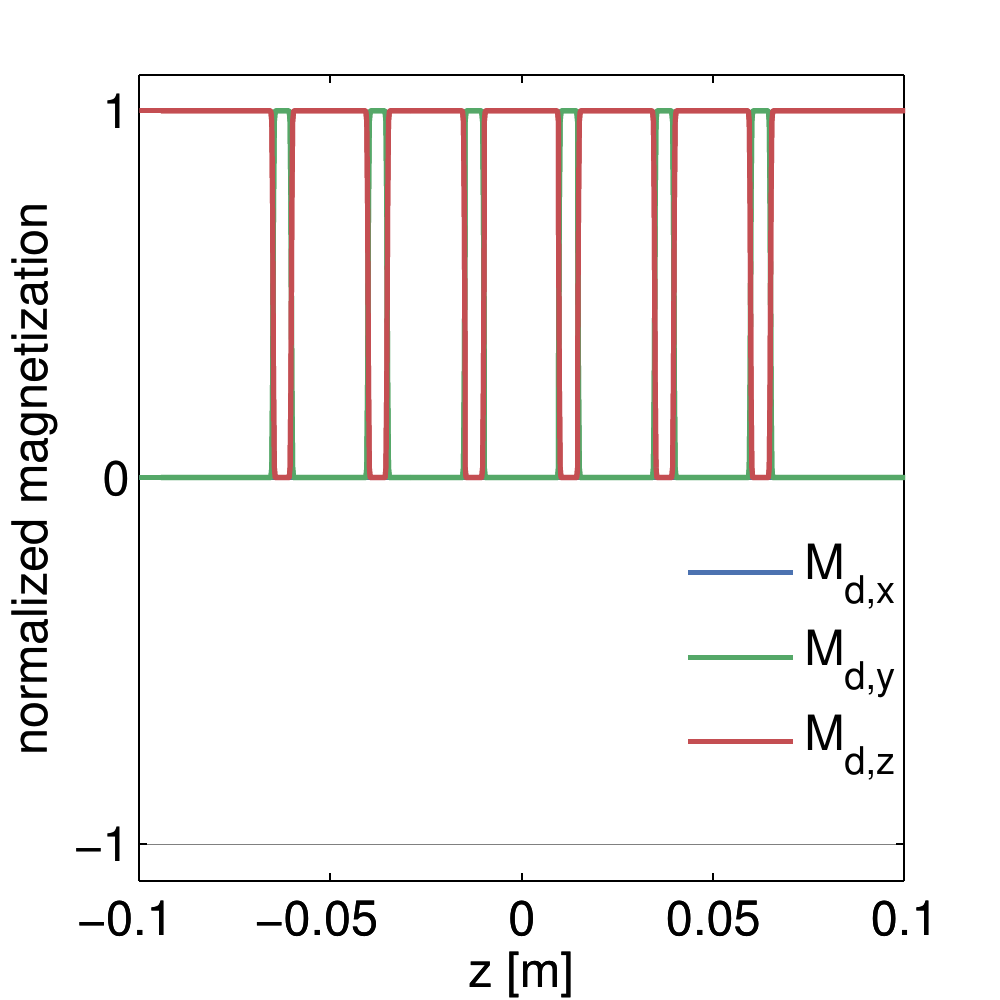}		
    }
    \hfill
    \subcaptionbox{${M_d}$ for phase shifted SMS 6 (zoom)\label{fig:initial:5_zoom}}{%
        \includegraphics[width=0.3\textwidth]{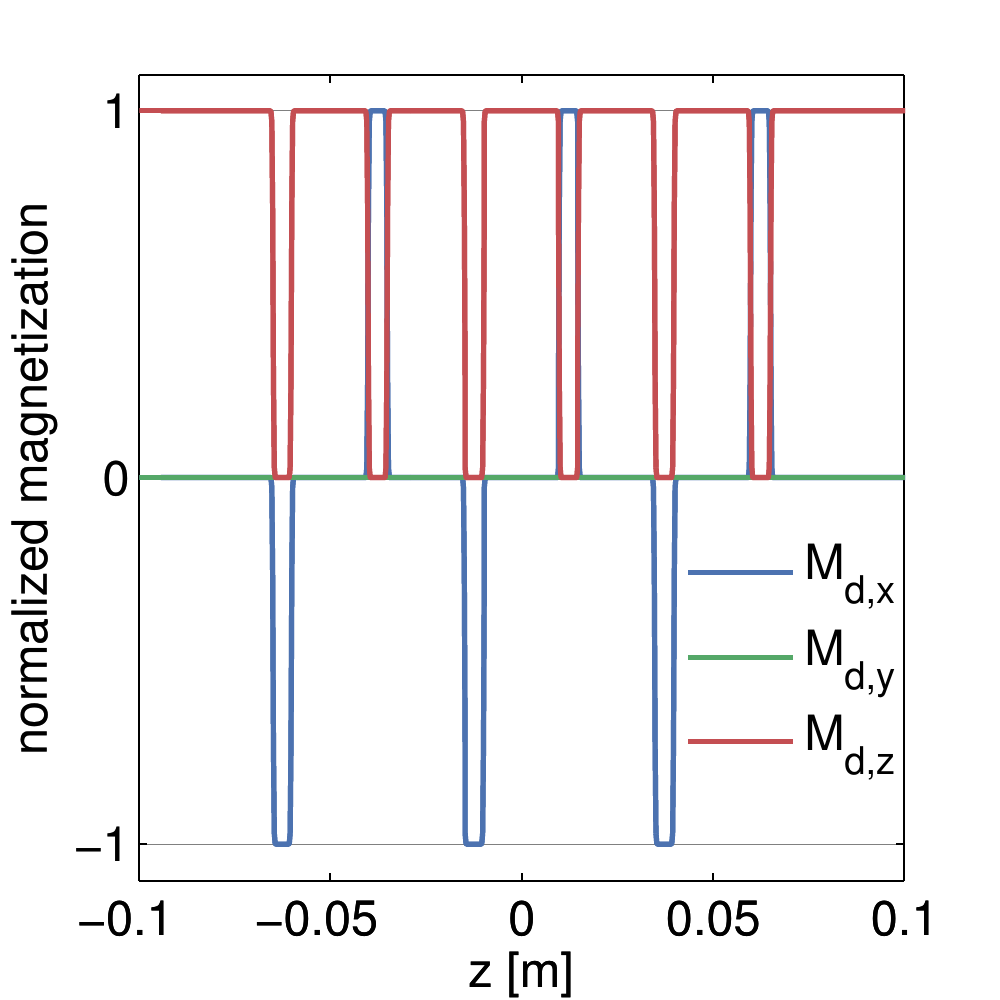}		
    }
    \hfill

    \caption{Desired magnetization for single- (a) and multi-slice (b,c) pulse design}\label{fig:initial}
\end{figure}

\paragraph{Single-slice excitation} 

To validate the design procedure, we compute an optimized pulse for a single slice of a given thickness $\Delta_w$ and a flip angle of \ang{90}, i.e., we set
\begin{equation}
    {\tilde M_d}(z) = \begin{cases} \left(0,\sin(\ang{90}),\cos(\ang{90})\right)^T &\text{if } |z|<\Delta_w/2,\\
        (0,0,1)^T &\text{else,}
    \end{cases}
\end{equation}
as visualized in \cref{fig:initial:3_zoom}.
To reduce Gibbs ringing, this vector is filtered before the optimization with a Gaussian kernel with a full width at half maximum of \SI{1.6}{\milli\metre}.
For comparison, an SLR pulse \cite{Pauly91,Lee07} with an identical temporal resolution and pulse duration is designed to the same specification (slice width, flip angle, full width at half maximum) using the Parks--McClellan (PM) algorithm \cite{Pauly91} with a \SI{1}{\percent} in- and out-of slice ripple as usual \cite{Pauly91} and a bandwidth of \SI{2.35}{\kilo\hertz}. 
To achieve a fully refocused magnetization, the refocusing area of the slice selective gradient for the SLR pulse is increased by \SI{4.1}{percent} compared to the OC pulse (\cref{fig:single:1}).

\paragraph{SMS excitation: phantom} RF pulses for the simultaneous excitation of two, four, and six equidistant rectangular slices with a flip angle of \ang{90} are computed, i.e., we set
\begin{equation}
{\tilde M_d}(z) = \begin{cases} \left(0,\sin(\ang{90}),\cos(\ang{90})\right)^T &\text{if $z$ in slice},\\
(0,0,1)^T &\text{if $z$ out of slice,}
    \end{cases}
\end{equation}
and apply Gauss filtering; see \cref{fig:initial:4_zoom} for the case of six slices.

Since PINS pulses are not suitable for axial or axial-oblique slice preference as they generate a periodic slice pattern extending outside the field of interest \cite{Ugurbil13}, the optimized pulses are compared with conventional SMS pulses obtained using superposed phase-shifted sinc-based excitation pulses, again for the same slice width, flip angle and full width at half maximum.

\paragraph{SMS excitation: in-vivo} 

Since multi-slice in-vivo imaging using slice-GRAPPA starts to suffer from g-factor problems for more than three slices, we modify the above-described SMS pulses using a CAIPIRINHA-based excitation pattern \cite{Breuer06}, which alternates two different pulses to achieve phase-shifted magnetization vectors in order to increase the spatial distance of aliased voxels.
Here, the first vector and pulse are identical to those designed for the phantom. For odd slice numbers, the second vector is modified by adding a phase term of $\pi$ to every second slice of the desired magnetization, i.e., 
\begin{equation}
{\tilde M_d}(z) = \begin{cases} \left(0,\pm\sin(\ang{90}),\cos(\ang{90})\right)^T &\text{if $z$ in odd/even slice},\\
(0,0,1)^T &\text{if $z$ out of slice}
    \end{cases}
\end{equation}
(before filtering). 
For even slice numbers, the transverse pattern has to be further shifted by $\frac{\pi}2$, i.e., 
\begin{equation}
{\tilde M_d}(z) = \begin{cases} \left(\pm\sin(\ang{90}),0,\cos(\ang{90})\right)^T &\text{if $z$ in even/odd slice},\\
(0,0,1)^T &\text{if $z$ out of slice,}
    \end{cases}
\end{equation}
see \cref{fig:initial:5_zoom} for the case of six slices. 
The additional phase shift is balanced before reconstruction by subtracting a phase of $\frac{\pi}{2}$ from every second phase-encoding line of the measured k-space data. 
Since typical relaxation times in the human brain are at least an order of magnitude bigger than the pulse duration, relaxation effects are neglected in the optimization. 

\bigskip

The starting point for the optimization is chosen in all cases as ${u}^0=[0,\dots,0]$.
The control cost parameter is fixed at $\alpha = 10^{-4}$ for both the single-slice and the multi-slice optimization.
The parameters in \cref{alg:trcgn} are set to $\mathrm{tol}_N = 10^{-9}$, $\mathrm{maxit}_N = 5$, $\mathrm{tol}_C = 10^{-6}$, $\mathrm{maxit}_{C} = 50$, $\rho_0 = 1$, $\rho_{\max} = 2$, $q=2$, $\sigma_1=0.03$, $\sigma_2 = 0.25$, $\sigma_3 = 0.7$. 
All calculations are performed on a workstation with a four-core \SI{64}{\bit} processor with \SI{3.1}{\giga\hertz} (Intel i\num{5}-\num{3350}P) and \SI{16}{\giga\byte} of RAM.

\subsection{Experimental validation}

Fully sampled experimental data for a phantom and a healthy volunteer were acquired on a \SI{3}{\tesla} MR scanner (Magnetom Skyra, Siemens Healthcare, Erlangen, Germany) using the built-in body coil to transmit the RF pulse. 
The MR signals were received using a body coil for the phantom experiments and a \num{32}-channel head coil for the in-vivo experiments.
A standard Cartesian GRE sequence was modified to import and apply external RF pulses.
By changing the read-out gradient from the frequency-axis to the slice direction, the excited slice can be measured and visualized. 
The single-slice excitation was measured using a water filled sphere with a diameter of \SI{170}{\milli\metre}. 
To acquire a high resolution in z-direction, we used a matrix size of \num{512x384} with a FOV of \SI{250x187}{\milli\metre} and a bandwidth of \SI{390}{\Hz}. The echo time was $T_E=$ \SI{5}{\milli\second} and the repetition time $T_R=$ \SI{2000}{\milli\second} to get fully relaxed magnetization before the next excitation.
The SMS phantom experiments were performed using a homogeneous cylinder phantom with diameter of \SI{140}{\milli\metre}, length of \SI{400}{\milli\metre}, and relaxation times $T_1=$ \SI{102}{\milli\second}, $T_2=$ \SI{81}{\milli\second}, and $T_2^*=$ \SI{70}{\milli\second}. 
The sequence parameters were $T_E=$ \SI{10}{\milli\second}, $T_R=$ \SI{1000}{\milli\second}, bandwidth \SI{390}{\Hz}, matrix size \num{512x288}, and a field of view of \SI{250x141}{\milli\metre}. 

To verify the in-vivo applicability, human brain images of a healthy volunteer were acquired using the above described GRE sequence modified to include the optimized CAIPIRINHA-based pulses. 
The sequence parameters were set to $T_E=$ \SI{10}{\milli\second}, $T_R=$ \SI{4000}{\milli\second}, bandwidth \SI{390}{\Hz}, matrix size \num{192x120} and FOV \SI{300x187}{\milli\metre}. 
After acquisition, the k-space data of the individual slices were separated using an offline slice-GRAPPA (\num{32} coils, kernel size of \num{4x4}) reconstruction \cite{Setsompop12,Cauley14}. The reference scans used in the slice-GRAPPA reconstruction were performed with the same sequence using an optimized single-slice pulse (not shown here). To decrease the scanning time, we acquired \num{25} k-space lines (\num[fraction-function=\sfrac]{1/5} of the full dataset) around the k-space center for each reference scan.
After this separation, a conventional Cartesian reconstruction was performed individually for each slice.

\section{Results}
\label{sec:res}

\paragraph{Single-slice excitation}

\begin{figure}[t]
    \centering
    \begin{minipage}[c]{0.66\textwidth}
        \captionsetup{format=hang,justification=centering}
        \centering
        \subcaptionbox{\mbox{optimized pulse $(u_x,u_y)$}, \mbox{gradient $G_z$ (fixed)}\label{fig:single:1}}{%
            \includegraphics[width=0.485\textwidth]{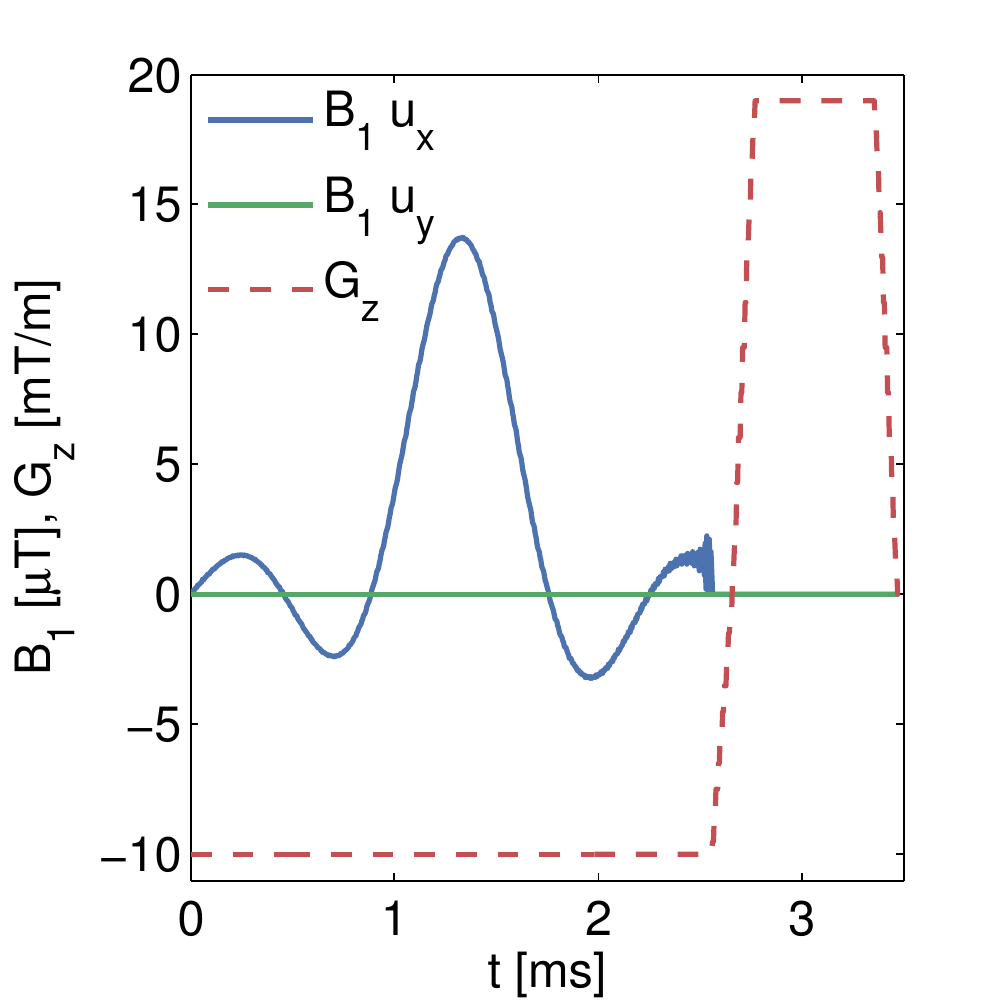}	
        }
        \hfill
        \subcaptionbox{\mbox{optimized transverse} \mbox{magnetization}\label{fig:single:2}}{%
            \includegraphics[width=0.485\textwidth]{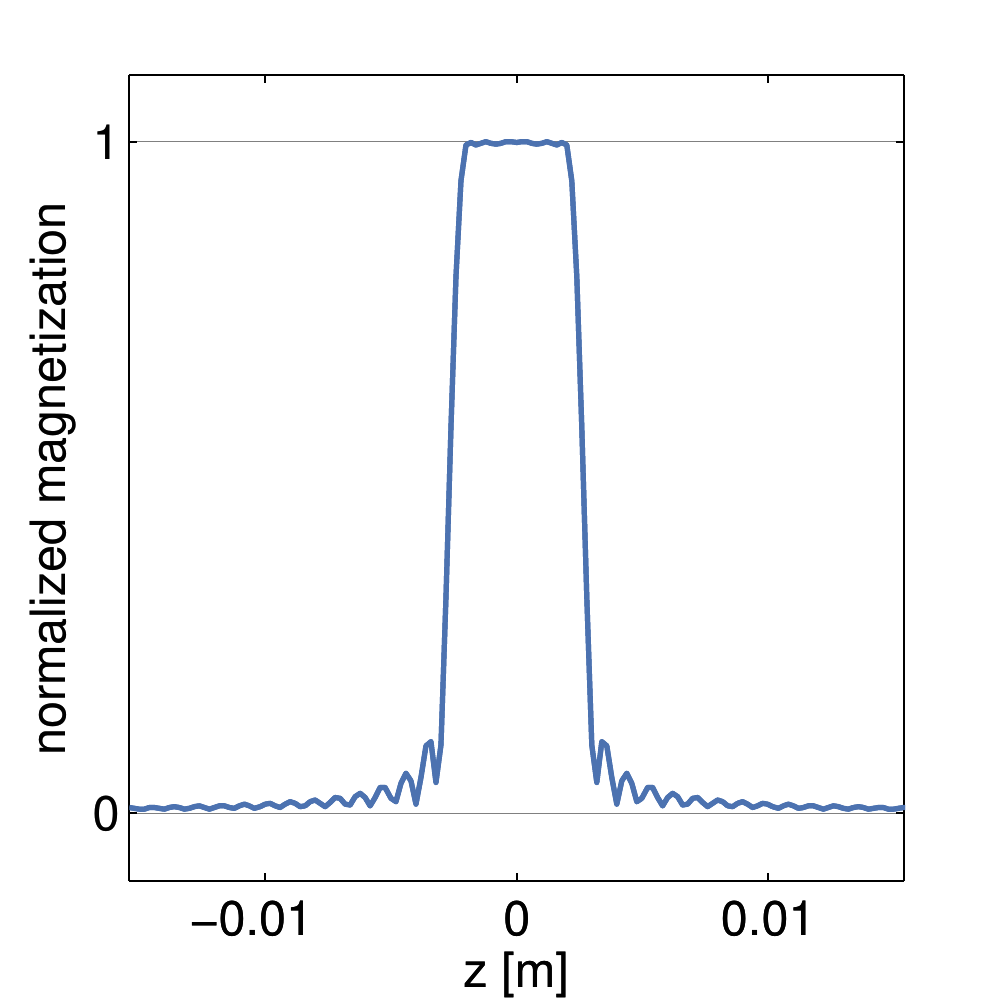}
        }    
        \\ 
        \hfill
        \subcaptionbox{reconstruction of excitation\label{fig:single:4}}{%
            \includegraphics[width=0.485\textwidth,clip,trim=0 30 10 30]{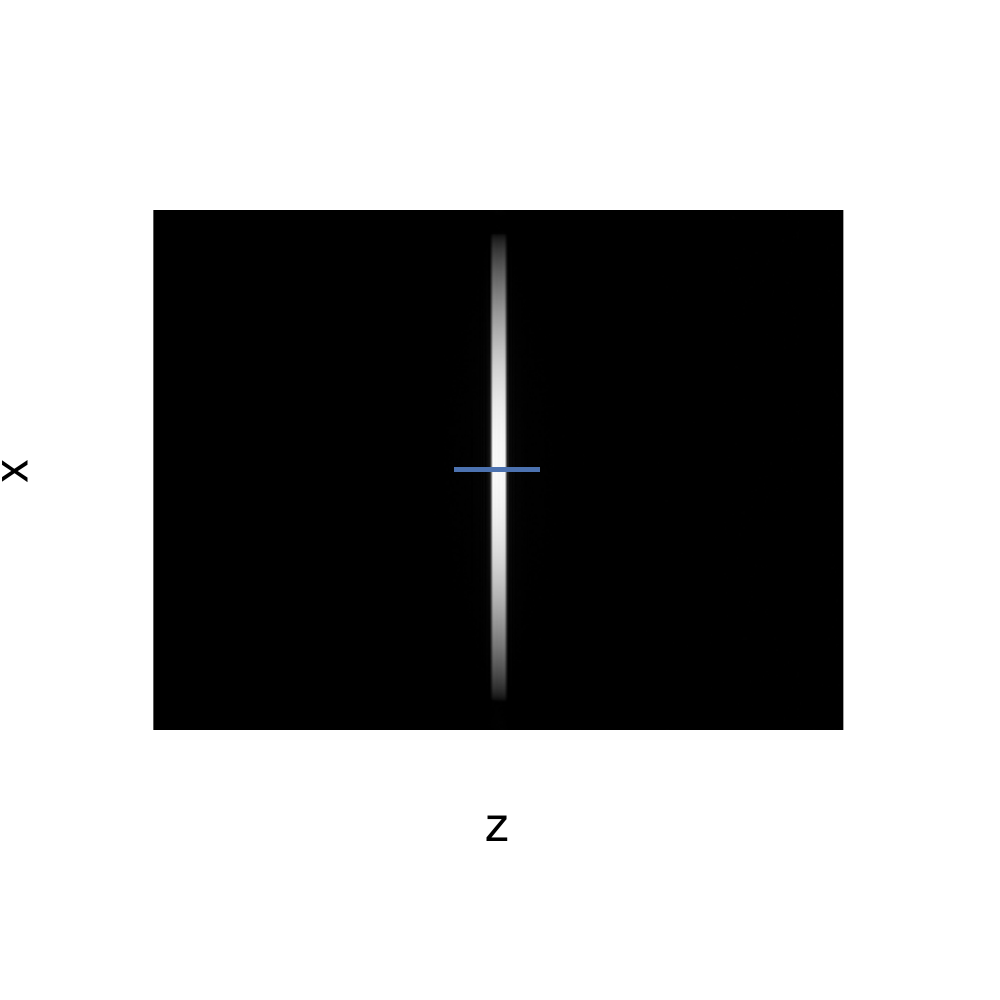}
        }
        \hfill
        \subcaptionbox{\mbox{measured transverse} magnetization\label{fig:single:3}}{%
            \includegraphics[width=0.485\textwidth]{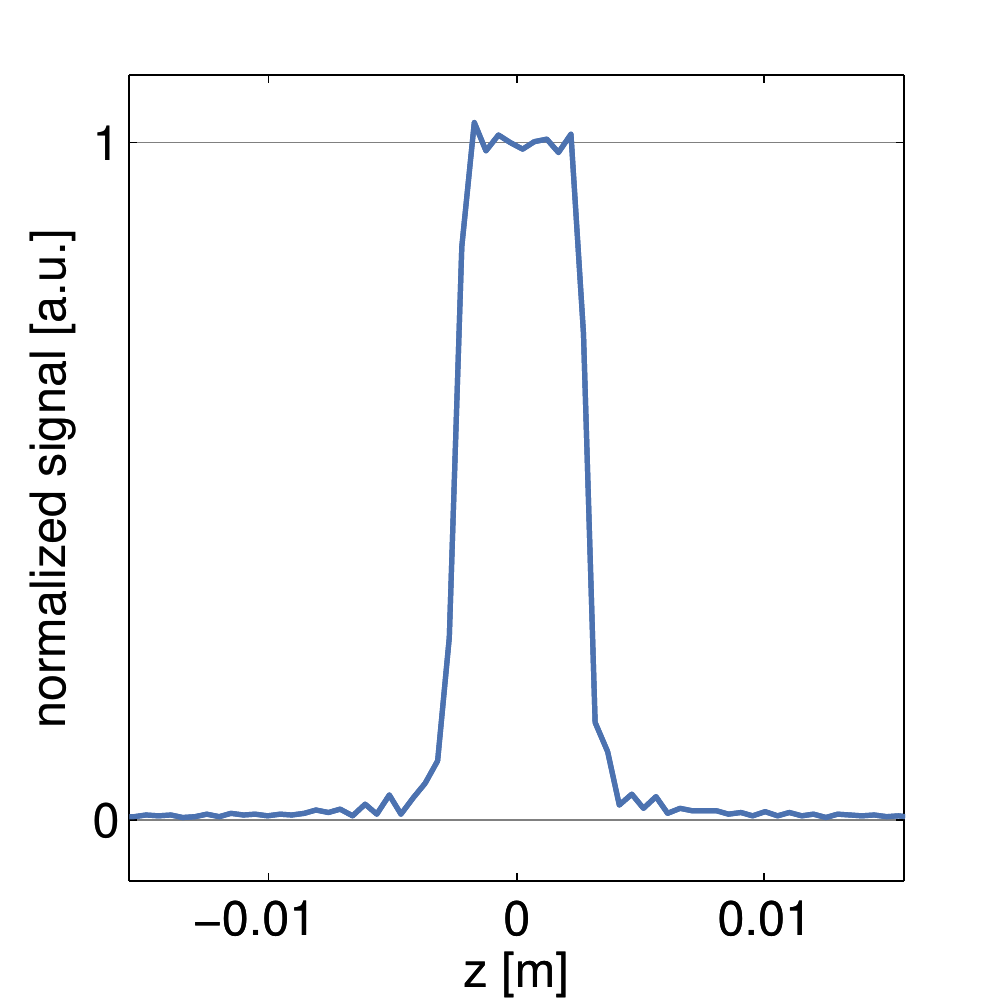}
        }
    \end{minipage}
    \caption{Optimized pulse and slice profile for single-slice excitation} \label{fig:single}
\end{figure}

\Cref{fig:single} shows the results of the design of an RF pulse for the excitation of a single slice of width $\Delta_w =$ \SI{5}{\milli\metre}; see \cref{fig:initial:3_zoom}. The computed pulse (after \num{4} Newton iterations and a total number of \num{28} CG steps taking \SI{989}{\second} on the above-mentioned workstation
is shown in \cref{fig:single:1}. (To indicate the sequence timing, the slice-selective gradient $G_z$ -- although not part of the optimization -- is shown dashed.) It can be seen that $u_x (t)$ is similar, but not identical, to a standard sinc shape, and that $u_y (t)$ is close to zero, which is expected due to the symmetry of the prescribed slice profile. \Cref{fig:single:2} contains a detail of the corresponding transverse magnetization ${M_{xy}}(T)=({M_x}(T)^2+{M_y}(T)^2)^{1/2}$ obtained from the numerical solution of the Bloch equation, which is confirmed by experimental phantom measurements in \cref{fig:single:4} and \cref{sub@fig:single:3}. 
Both simulation and measurement show an excitation with a steep transition between the in- and out-of-slice regions and a homogeneous flip angle distribution across the target slice. 

\begin{figure}[t]
    \centering
    \begin{minipage}[c]{0.66\textwidth}
        \centering
        \subcaptionbox{simulated magnetization (OC)\label{fig:comp:5}}{%
            \includegraphics[width=0.485\textwidth]{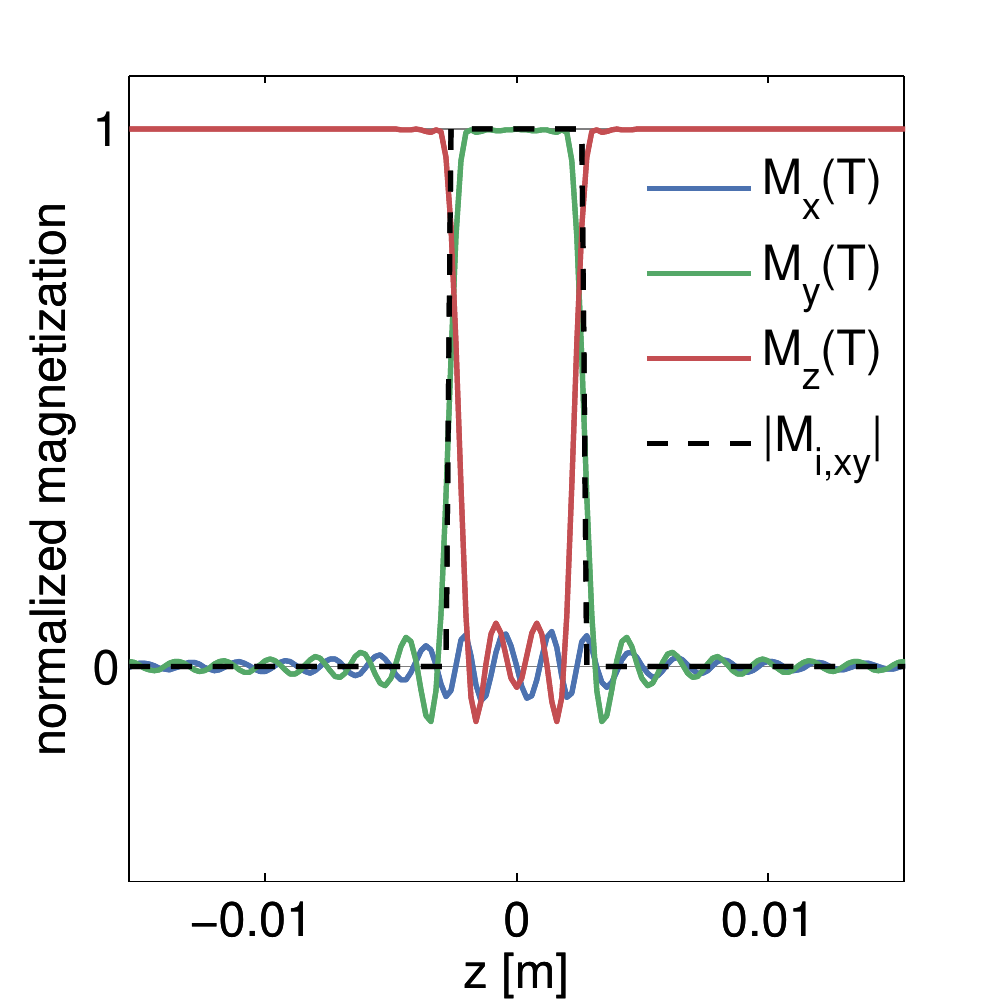}		
        }
        \hfill
        \subcaptionbox{simulated magnetization (SLR)\label{fig:comp:6}}{%
            \includegraphics[width=0.485\textwidth]{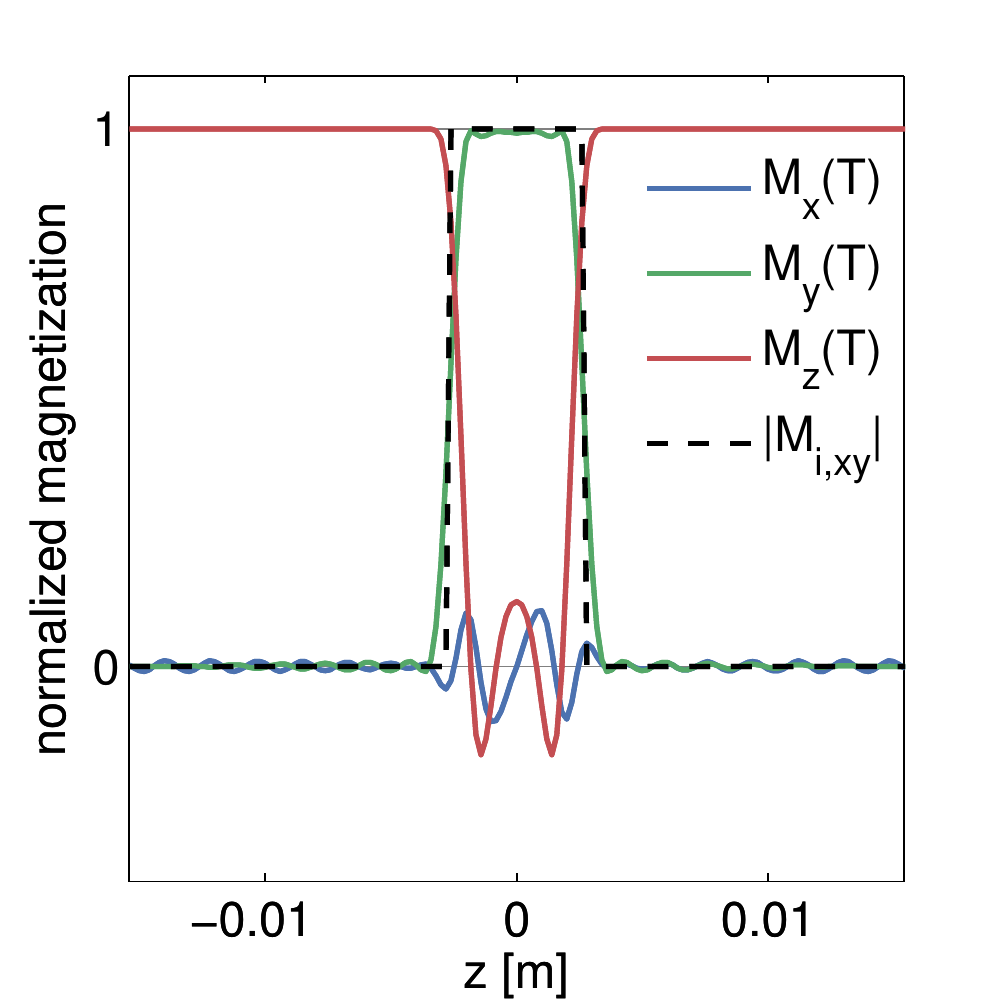}  
        }\\
        \hfill

        \subcaptionbox{in-slice phase\label{fig:comp:3}}{%
            \includegraphics[width=0.485\textwidth]{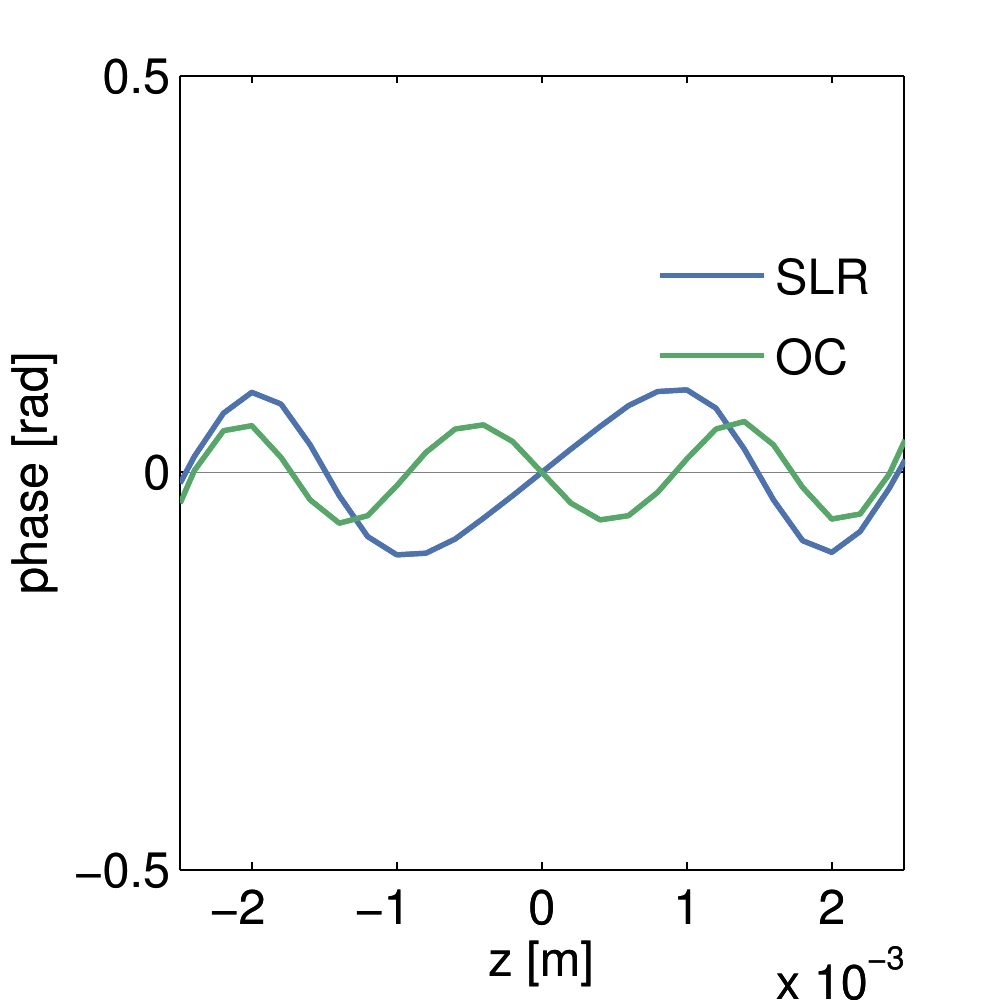}
        }
        \hfill
        \subcaptionbox{ideal magnetization\label{fig:comp:4}}{%
            \includegraphics[width=0.485\textwidth]{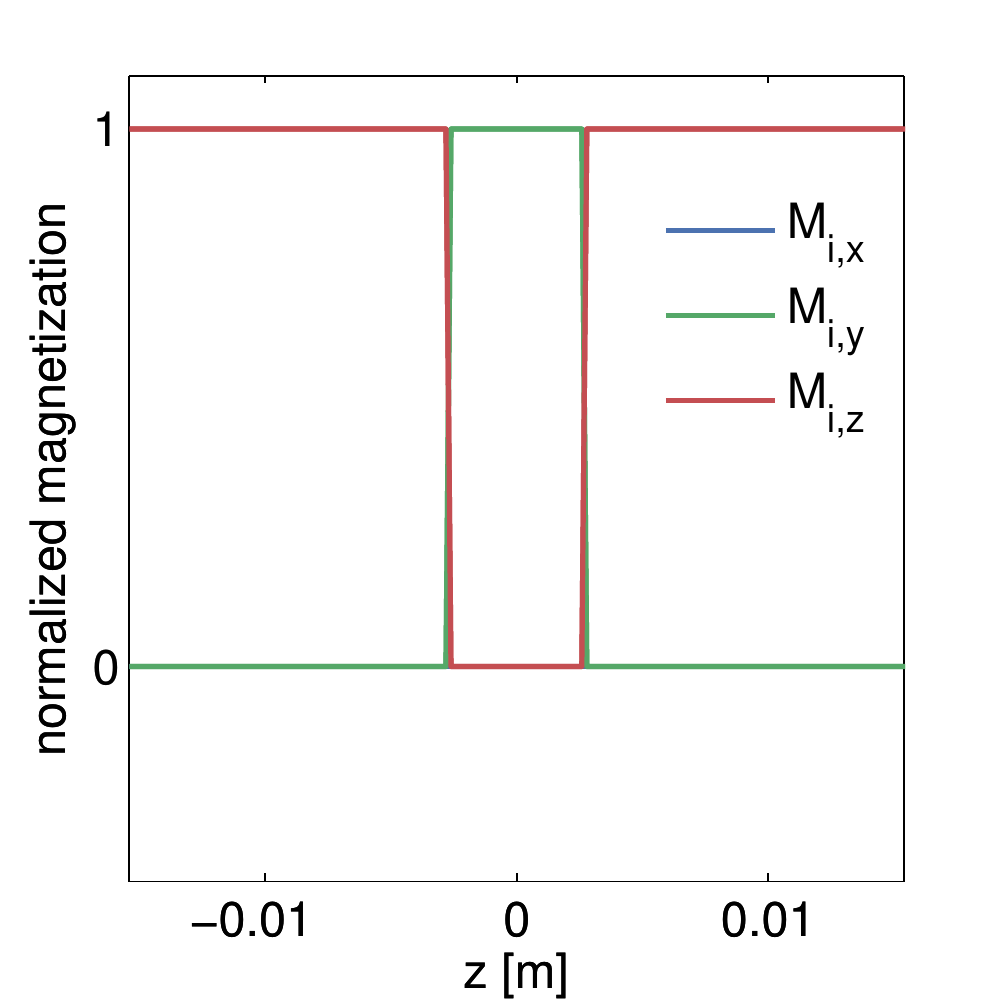}
        }
    \end{minipage}
    \caption{Comparison of SLR and OC pulse}\label{fig:comp}
\end{figure}

\Cref{fig:comp} compares the optimized (OC) pulse with a standard SLR pulse by showing details of the corresponding simulated magnetizations (\cref{fig:comp:5} for OC and \cref{fig:comp:6} for SLR; in both cases the targeted ideal magnetization is shown dashed). 
It can be seen that the in-slice magnetization of the optimized pulse has oscillations of higher frequency but of much smaller amplitude than that of the SLR pulse. This becomes especially visible when comparing the resulting in-slice phases (\cref{fig:comp:3}).

This is achieved by allowing higher ripples close to the slice while decreasing the amplitude monotonically away from the slice. (Note that only a small central segment of this region is shown in the figures.) This leads to the total root mean squared error (RMSE) and the mean absolute error (MAE) with the ideal rectangular magnetization pattern (\cref{fig:comp:4}) matching the full width at half maximum of both pulses being smaller for the OC pulse (\num{1.46e-2} and \num{1.10e-4}, respectively) compared to the SLR pulse (\num{1.62e-2} and \num{2.27e-4}, with an equal power demand for both pulses. 

\paragraph{SMS excitation: phantom}

\enlargethispage*{1cm}

\begin{figure}[p]
    \centering
    \subcaptionbox{optimized pulse (2 slices)\label{fig:sms:a2}}{%
        \includegraphics[width=0.3\textwidth]{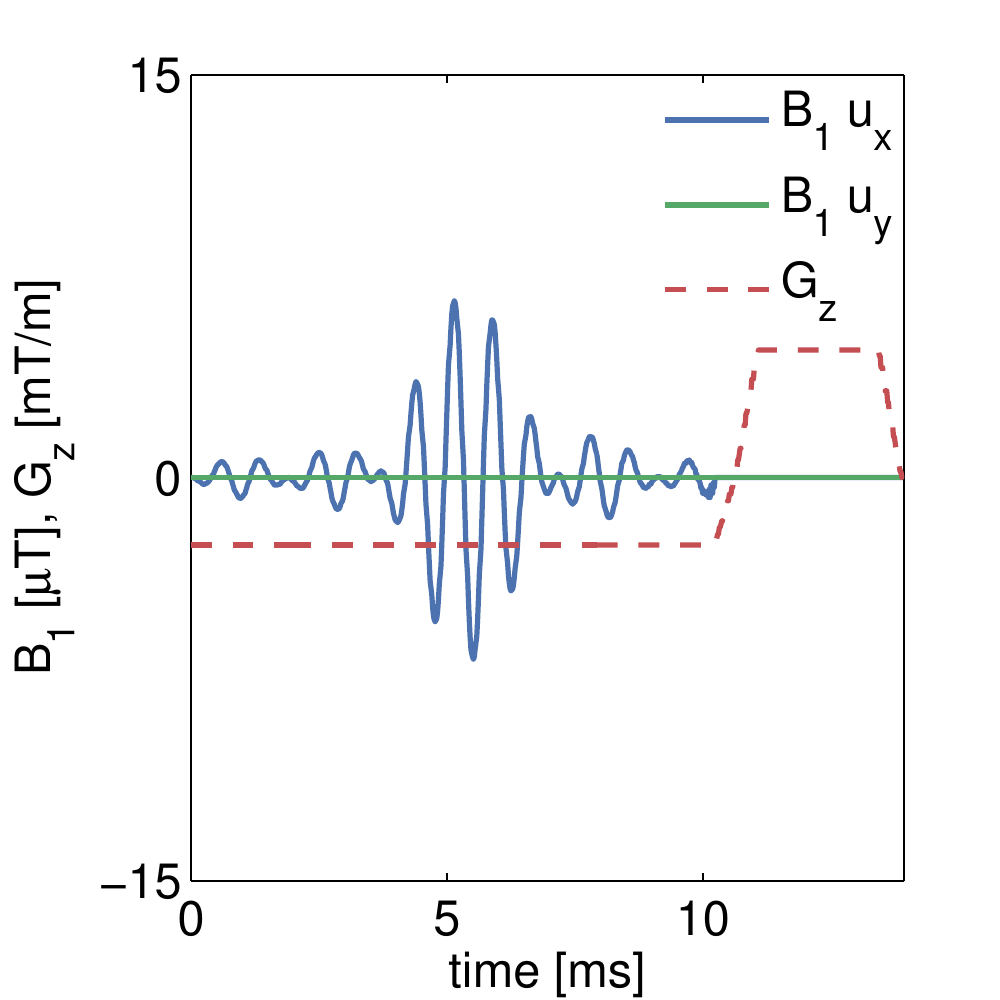}
    }
    \hfill
    \subcaptionbox{optimized pulse (4 slices)\label{fig:sms:a4}}{%
        \includegraphics[width=0.3\textwidth]{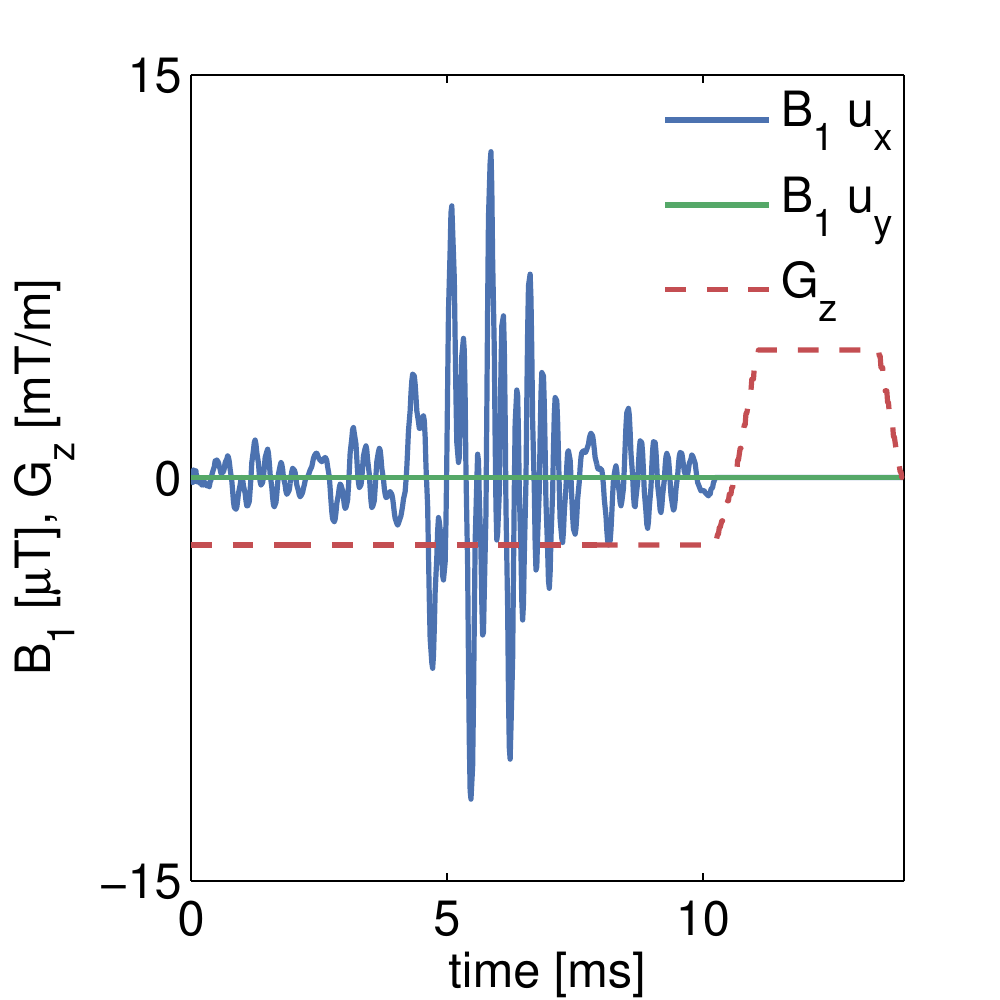}
    }
    \hfill
    \subcaptionbox{optimized pulse (6 slices)\label{fig:sms:a6}}{%
        \includegraphics[width=0.3\textwidth]{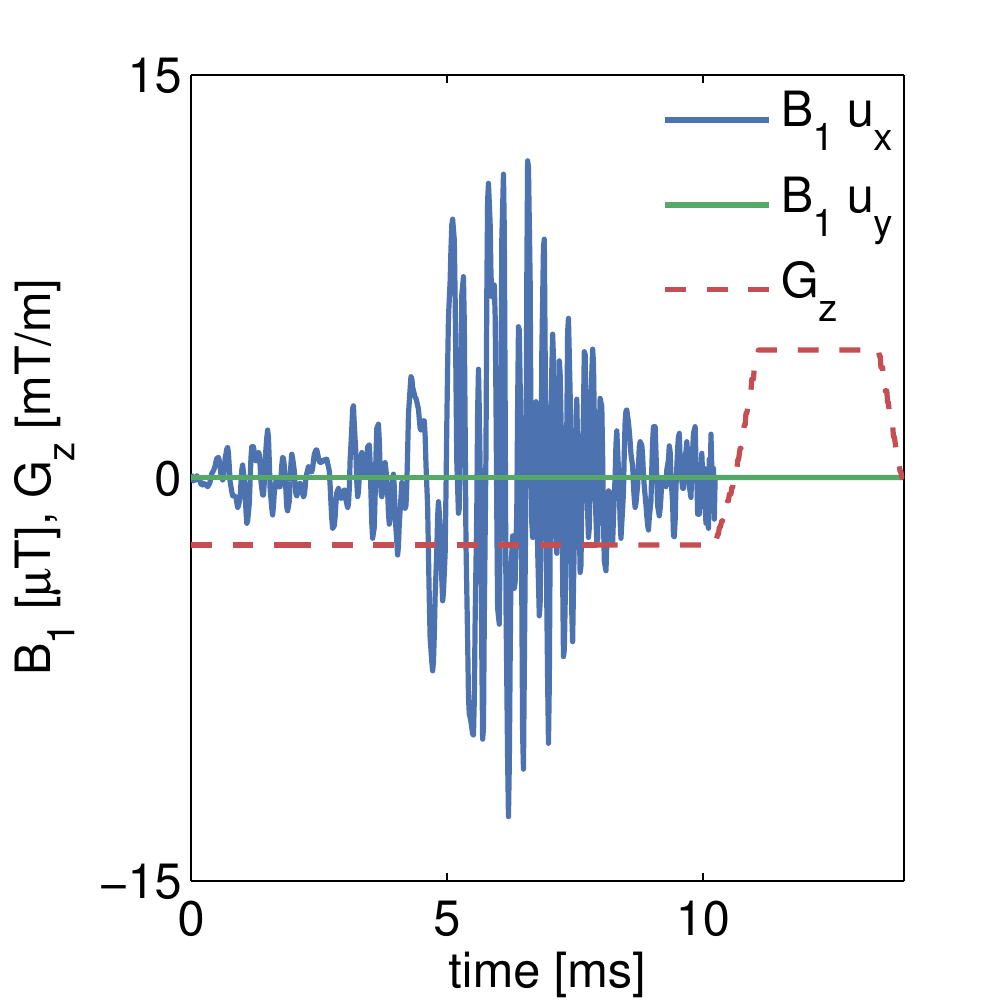} 
    }

    \subcaptionbox{optimized slice profile (2)\label{fig:sms:b2}}{%
        \includegraphics[width=0.3\textwidth]{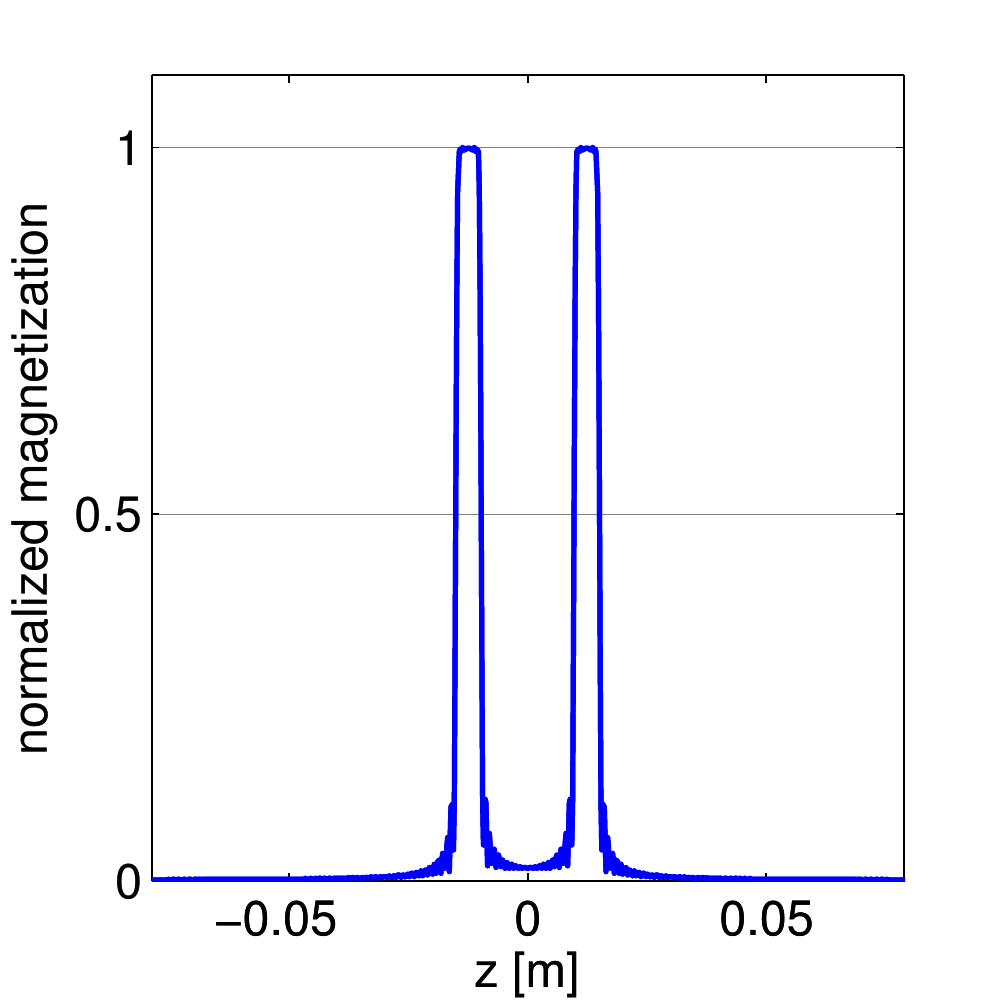}  
    }
    \hfill
    \subcaptionbox{optimized slice profile (4)\label{fig:sms:b4}}{%
        \includegraphics[width=0.3\textwidth]{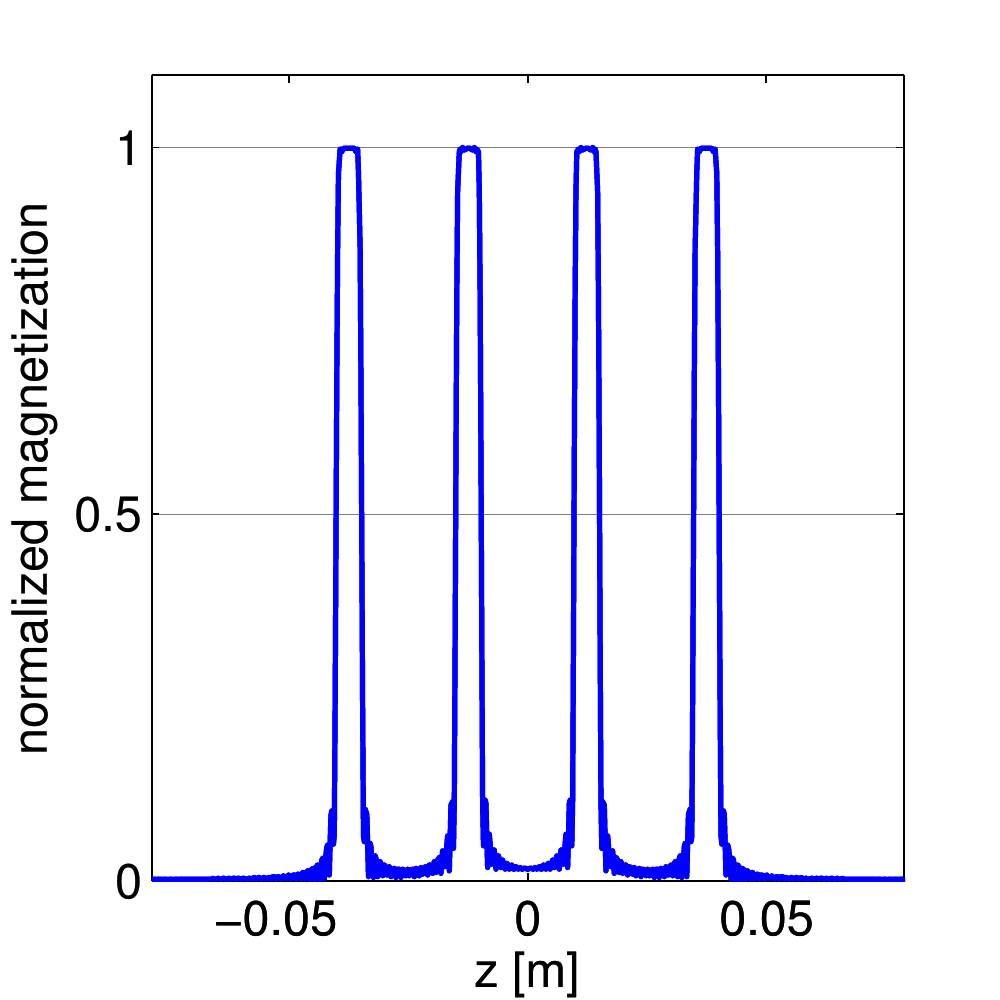}  
    }
    \hfill
    \subcaptionbox{optimized slice profile (6)\label{fig:sms:b6}}{%
        \includegraphics[width=0.3\textwidth]{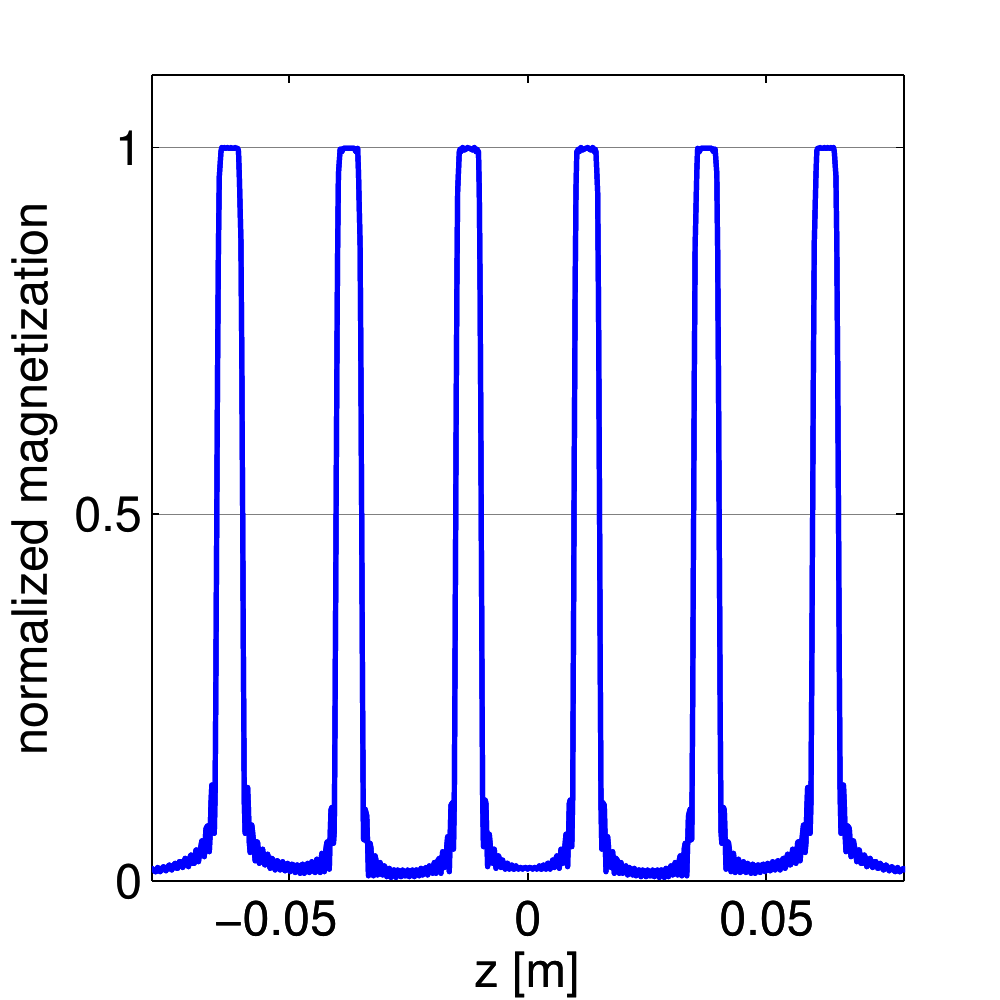} 
    }

    \subcaptionbox{reconstructed excitation (2)\label{fig:sms:c2}}{%
        \includegraphics[width=0.3\textwidth,clip,trim=0 50 30 50]{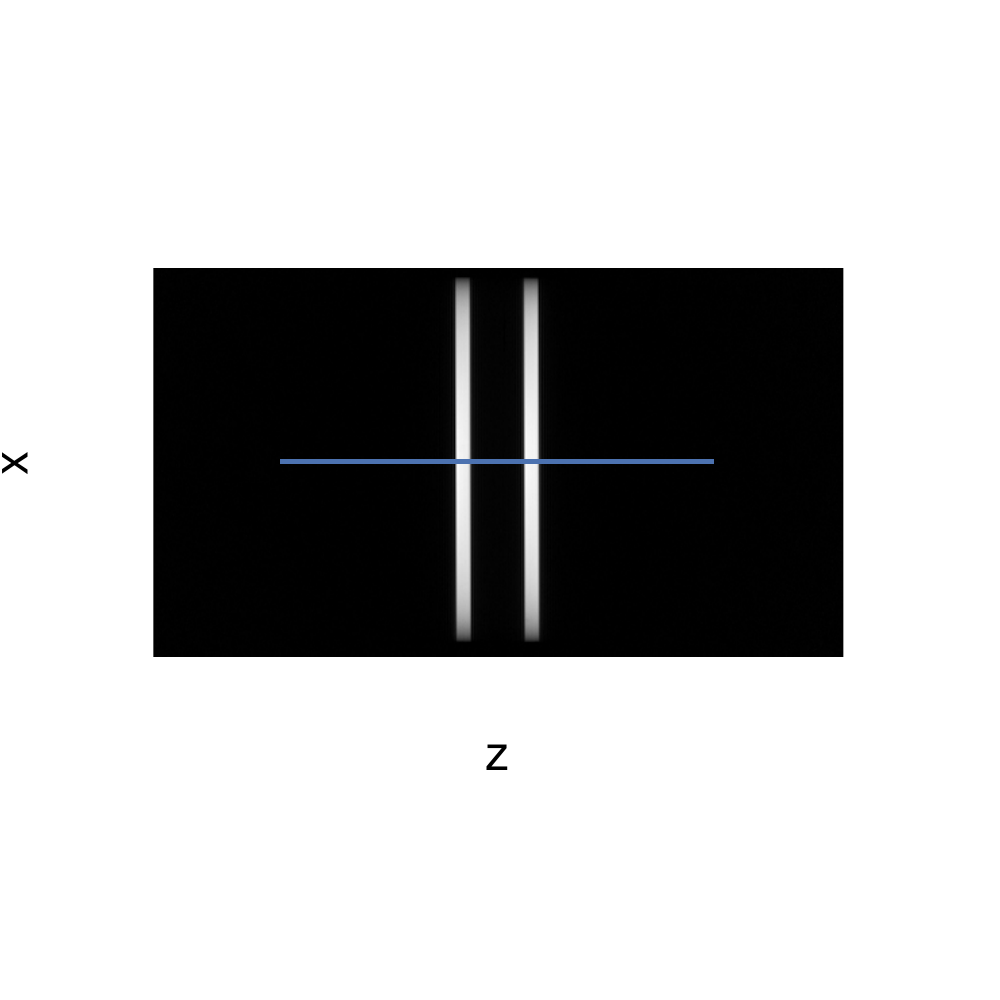}  
    }
    \hfill
    \subcaptionbox{reconstructed excitation (4)\label{fig:sms:c4}}{%
        \includegraphics[width=0.3\textwidth,clip,trim=0 50 30 50]{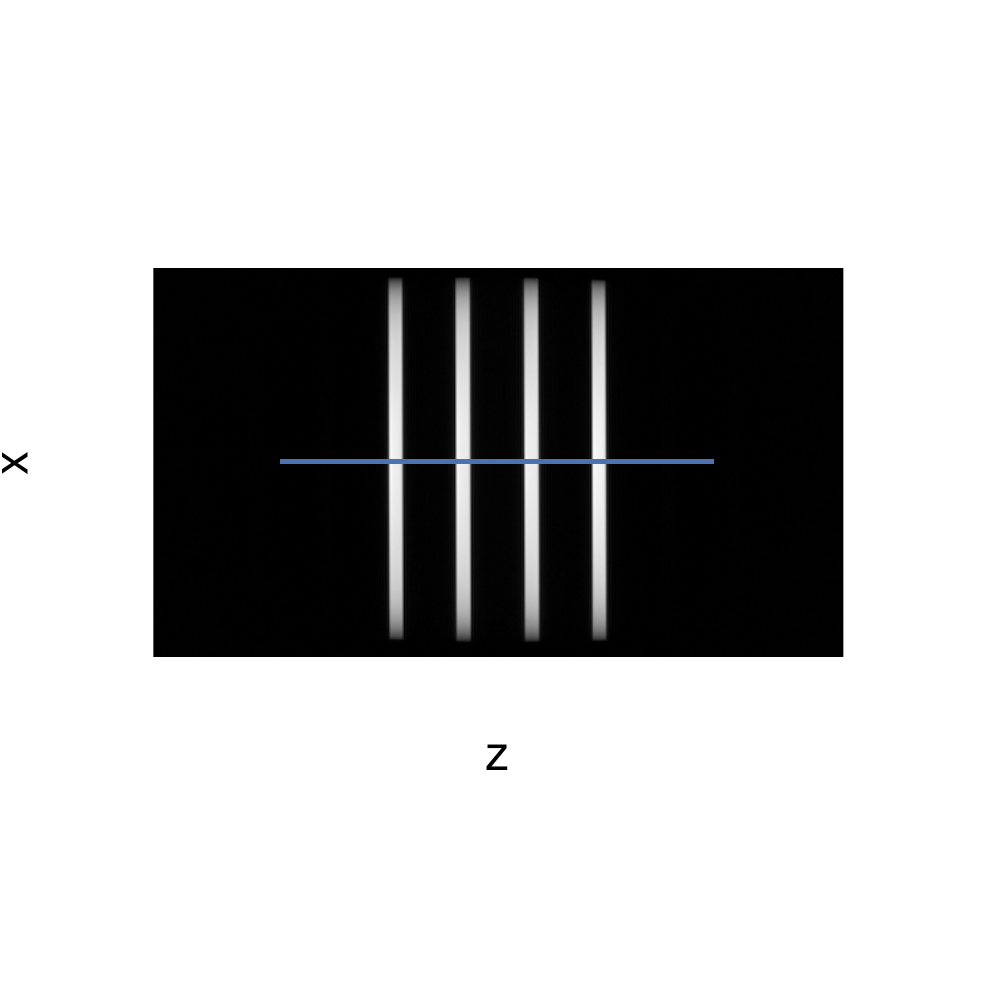} 
    }
    \hfill
    \subcaptionbox{reconstructed excitation (6)\label{fig:sms:c6}}{%
        \includegraphics[width=0.3\textwidth,clip,trim=0 50 30 50]{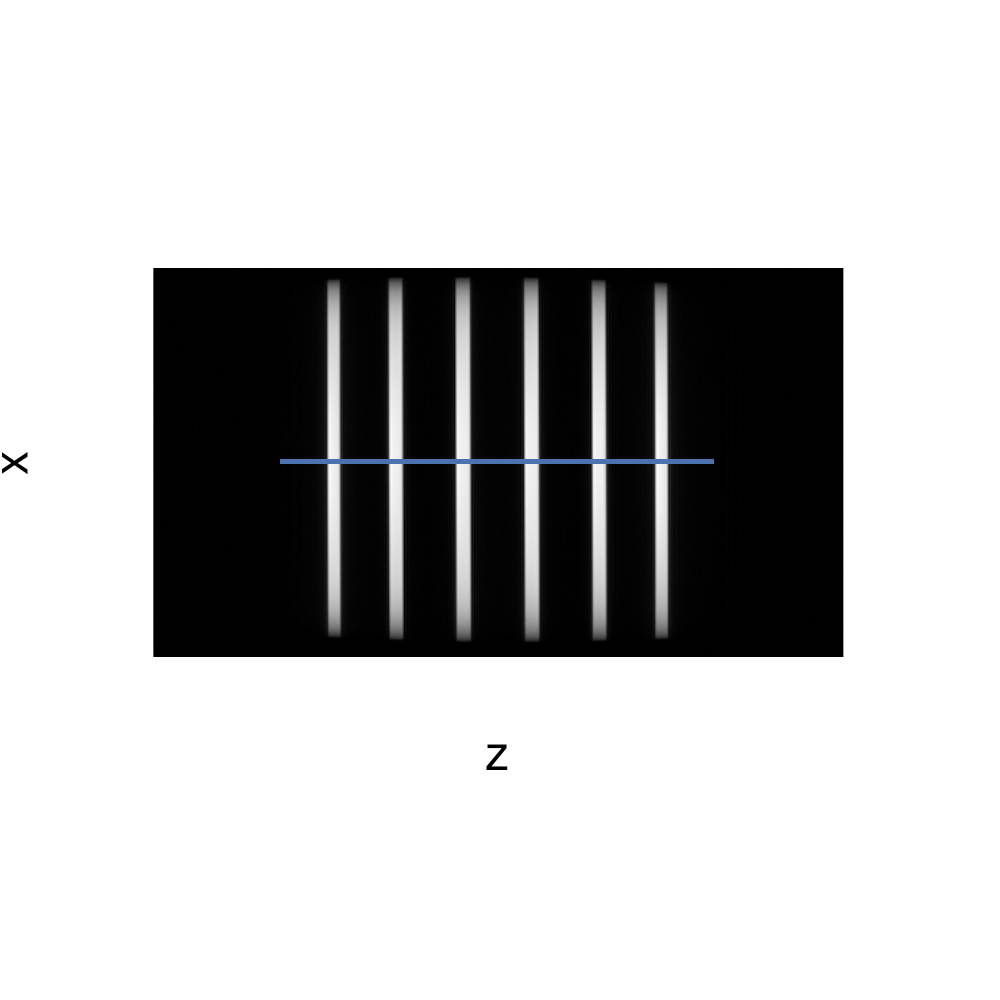} 
    }

    \subcaptionbox{measured slice profile (2)\label{fig:sms:d2}}{%
        \includegraphics[width=0.3\textwidth]{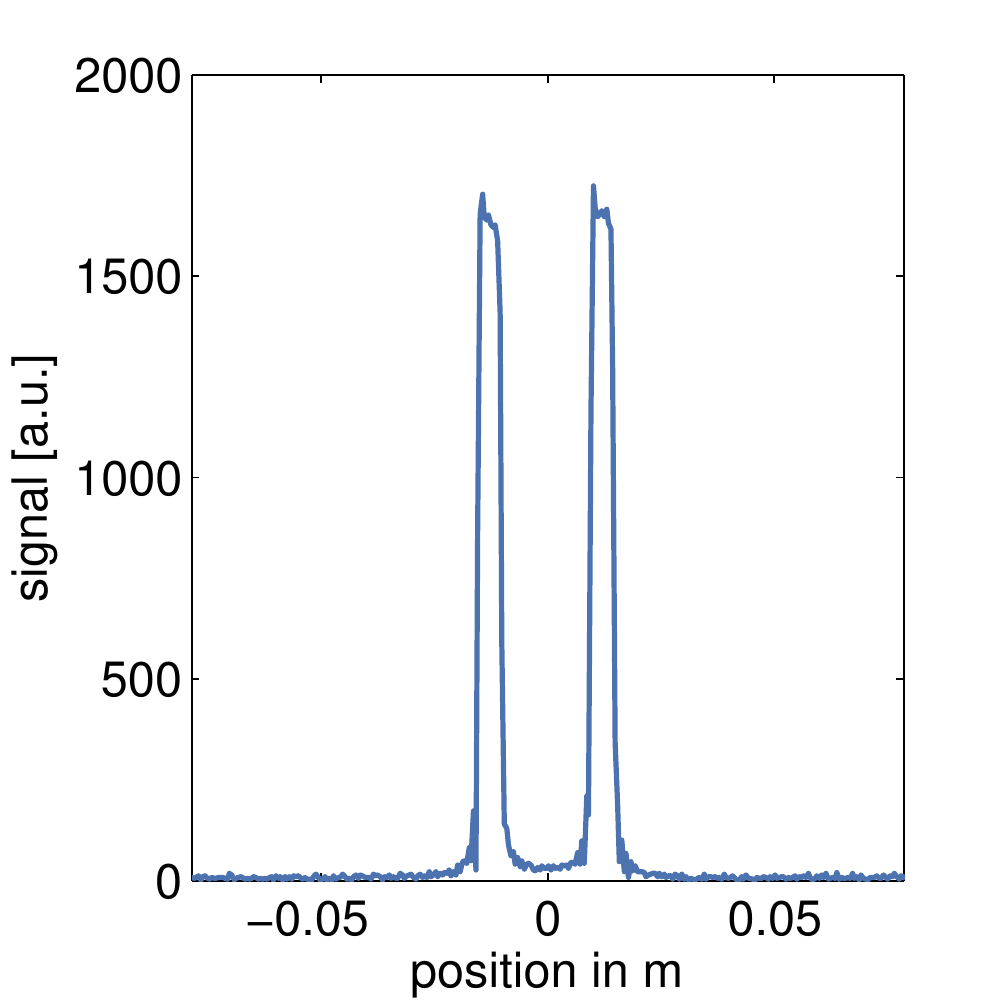}
    }
    \hfill
    \subcaptionbox{measured slice profile (4)\label{fig:sms:d4}}{%
        \includegraphics[width=0.3\textwidth]{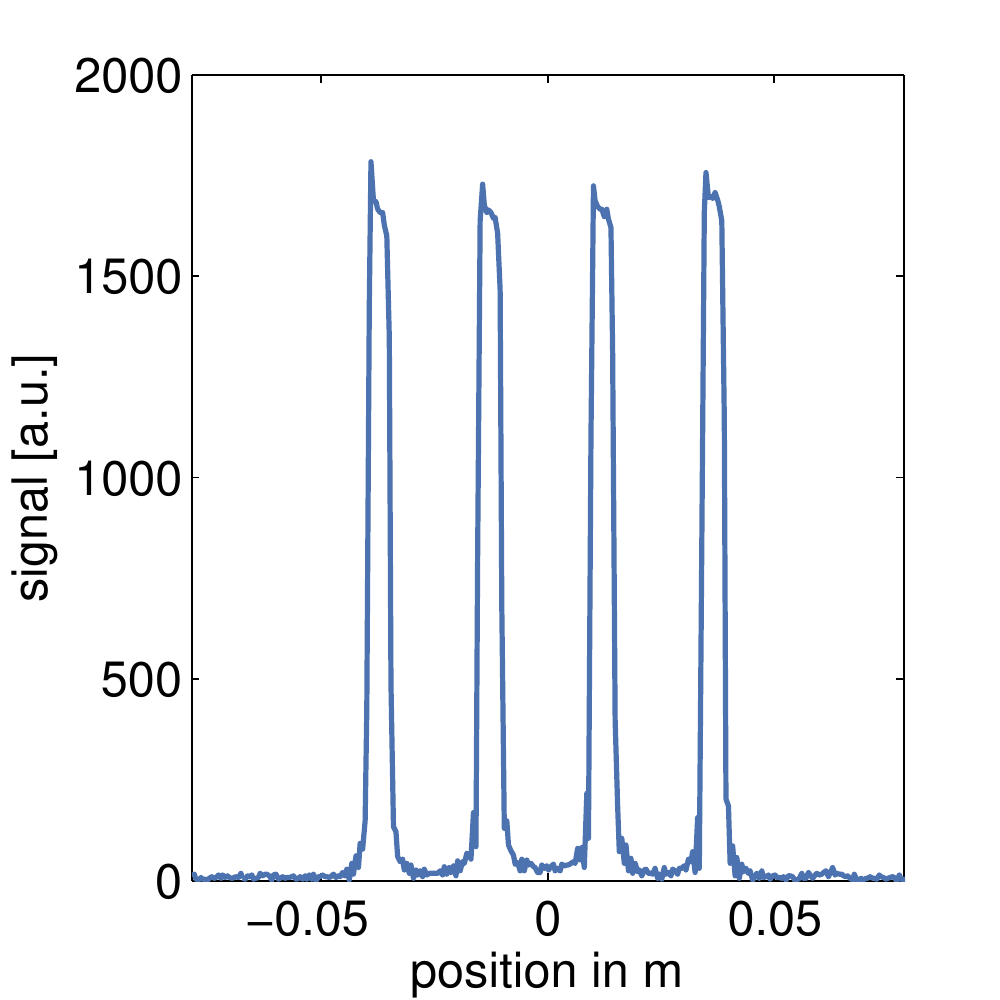} 
    }
    \hfill
    \subcaptionbox{measured slice profile (6)\label{fig:sms:d6}}{%
        \includegraphics[width=0.3\textwidth]{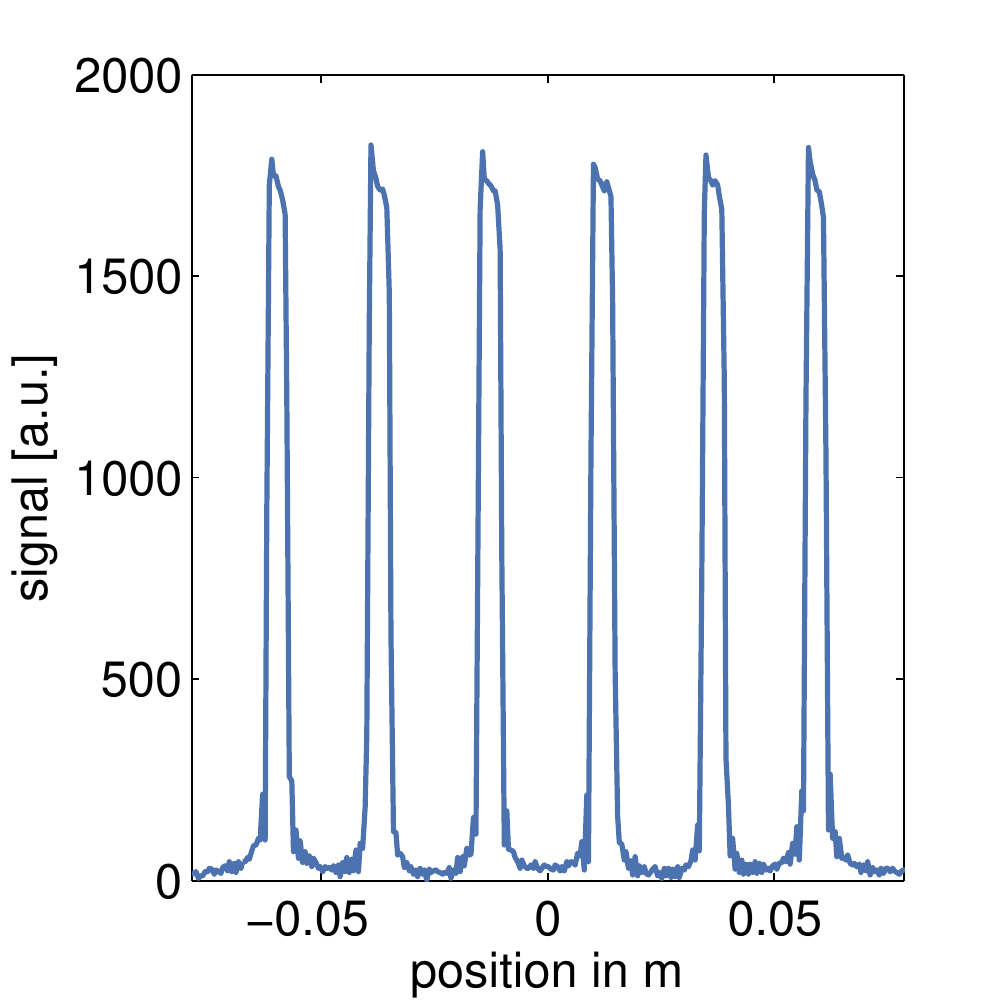}
    }
    \caption{Optimized pulses and slice profiles for SMS excitation (phantom)}\label{fig:sms}
\end{figure}

\Cref{fig:sms} shows the results of the design of RF pulses for simultaneous excitation of two, four and six equidistant slices with a separation of \SI{25}{\milli\metre} and a thickness $\Delta_w=$ \SI{5}{\milli\metre}; see \cref{fig:initial:4_zoom}.
The computational effort in all cases is similar to that in the single-slice case. 
The corresponding computed pulses are shown in \crefrange{fig:sms:a2}{sub@fig:sms:a6}. A graphical analysis shows that instead of higher amplitudes, the optimization distributes the total RF power (which increases with the number of slices) more uniformly over the pulse length. 
A central section of the corresponding optimized slice profiles are given in \crefrange{fig:sms:b2}{sub@fig:sms:b6}. It can be seen that all slices have a sharp profile which does not deteriorate as the number of slices increases (although it decreases slightly farther from the center and the bandwidth is slightly reduced). 
These results are validated by the experimental phantom measurements using the computed pulses: \crefrange{fig:sms:c2}{sub@fig:sms:c6} show the reconstructed excitation inside the phantom, while 
\crefrange{fig:sms:d2}{sub@fig:sms:d6} show the measured slice profiles along a cut parallel to the $x$-axis in the center of the previous images.

A quantitative comparison of SLR and OC-based SMS pulses from one to six simultaneous slices is given in \cref{tab:exp2}, which shows both the power requirement of the computed pulses, both in total B$_1$ energy
\begin{equation}
    \|{B_{1,x}}\|_2^2 = \int_0^T |B_1 u_x(t)|^2\,dt    
\end{equation}
and in peak B$_1$ amplitude
\begin{equation}
    \|{B_{1,x}}\|_\infty = \max_{t\in[0,T]} |B_1 u_x(t)|,
\end{equation}
as well as the mean absolute error (MAE) with respect to the ideal (unfiltered) slice profiles for the in-slice and the out-of-slice regions.
While both methods lead to a linear increase of the total energy with the number of slices, the peak amplitude increases more slowly for the OC pulses than for the conventional pulses.
Furthermore, we remark that the peak B$_1$ amplitude for four, five and six slices remain similar.
Regarding the corresponding slice profiles, the OC pulses lead to a significantly lower MAE in both the in-slice and out-of-slice regions compared to the SLR pulses. Visual inspection of \crefrange{fig:sms:b2}{sub@fig:sms:b6} shows that this is due to the fact that the out-of-slice ripples are concentrated around the in-slice regions while quickly decaying away from them.
\begin{table}[t]
    \caption{Comparison of B$_1$ power and the mean absolute error (MAE) of the transverse magnetization after excitation for conventional and OC based SMS pulses}\label{tab:exp2}
    \centering
    \begin{tabular}{S[table-format=1.0]
            S[table-format=2.1]
            S[table-format=2.1]
            S[table-format=2.1]
            S[table-format=2.2]
            S[table-format=1.3]
            S[table-format=1.3]
            S[table-format=1.4]
        S[table-format=1.4]}
        \toprule
        & \multicolumn{2}{c}{$\|B_{1,x}\|_2^2$ } &  \multicolumn{2}{c}{$\|B_{1,x}\|_{\infty}$}    
        & \multicolumn{2}{c}{MAE in-slice}       &   \multicolumn{2}{c}{MAE out-of-slice}\\
        & \multicolumn{2}{c}{[a.u.]} &  \multicolumn{2}{c}{[\si{\micro\tesla}]}
        & \multicolumn{2}{c}{[a.u.]} &  \multicolumn{2}{c}{[a.u.]}\\
        \cmidrule(lr){2-3} 
        \cmidrule(lr){4-5}
        \cmidrule(lr){6-7}
        \cmidrule(lr){8-9}

        {slices} & {conv} & {OC}  & {conv} & {OC}  & {conv} & {OC}  & {conv} & {OC}  \\
        \midrule
        1        & 19.5   & 19.5  & 3.5    & 3.49  & 0.062  & 0.052 & 0.0039 & 0.0014\\
        2        & 38.9   & 38.1  & 7.0    & 6.78  & 0.060  & 0.052 & 0.0040 & 0.0018\\
        3        & 58.4   & 57.2  & 10.5   & 10.02 & 0.054  & 0.053 & 0.0039 & 0.0030\\
        4        & 77.9   & 76.3  & 14.0   & 12.13 & 0.065  & 0.045 & 0.0086 & 0.0031\\
        5        & 97.3   & 95.5  & 17.5   & 11.38 & 0.059  & 0.053 & 0.0078 & 0.0051\\
        6        & 116.8  & 113.9 & 21.0   & 12.63 & 0.068  & 0.053 & 0.0075 & 0.0067\\
        \bottomrule
    \end{tabular}
\end{table}

Finally, we illustrate the influence of the regularization parameter $\alpha$ in \cref{tab:regalpha}, where the root of mean square error (RMSE), the total B$_1$ energy as well as the B$_1$ peak of the OC SMS 6 pulses is shown for different values of the control cost parameter $\alpha$. As can be seen, a bigger $\alpha$ leads to an increase in the error between desired and controlled magnetization while both the total B$_1$ power and the peak B$_1$ amplitude are reduced, although these effects amount to less than \SI{20}{percent} over a range of parameters spanning two orders of magnitude. This demonstrates that the results presented here are robust with respect to the choice of the control cost parameter.
\begin{table}[t]
    \caption{Comparison of RMSE, B$_1$ power and B$_1$ peak for different values of $\alpha$}\label{tab:regalpha}
    \centering
    \begin{tabular}{S[table-format=1.0e+1] 
            S[table-format=1.3e+1] 
            S[table-format=3.1]
            S[table-format=2.2]
        }
        \toprule
        {$\alpha$} & {RMSE}   & {$\|B_{1,x}\|_2^2$} & {$\|B_{1,x}\|_{\infty}$} \\
        {[a.u.]}   & {[a.u.]} & {[a.u.]}            & {[\si{\micro\tesla}]}     \\
        \midrule
        1e-5       & 2.374e-2 & 117.0               & 12.75   \\
        5e-5       & 2.375e-2 & 115.1               & 12.71   \\
        1e-4       & 2.377e-2 & 113.9               & 12.62   \\
        5e-4       & 2.437e-2 & 106.7               & 12.14   \\
        1e-3       & 2.591e-2 & 98.9                & 11.63   \\
        \bottomrule
    \end{tabular}
\end{table}

\paragraph{SMS excitation: in-vivo}

\begin{figure}[t]
    \centering
    \begin{minipage}[c]{0.7\textwidth}
        \subcaptionbox{optimized pulse (no shift)  \label{fig:method1:1}}{%
            \includegraphics[width=0.44\textwidth]{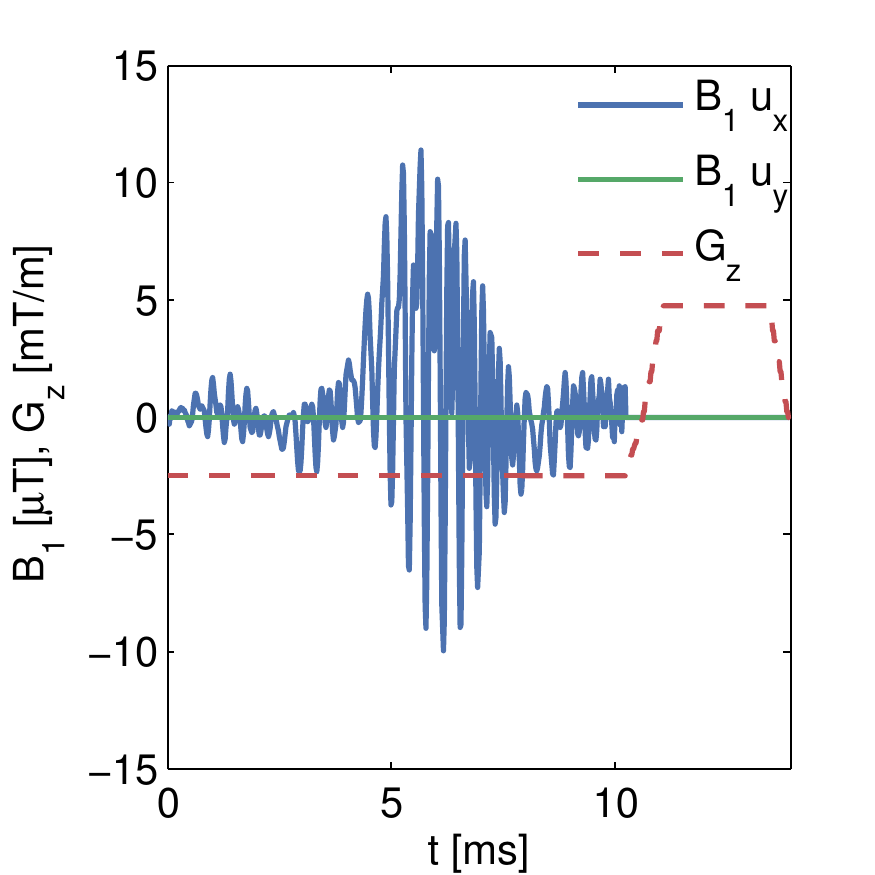}
        }
        \hfill
        \subcaptionbox{optimized slice profile of (a) \label{fig:method1:2}}{%
            \includegraphics[width=0.44\textwidth]{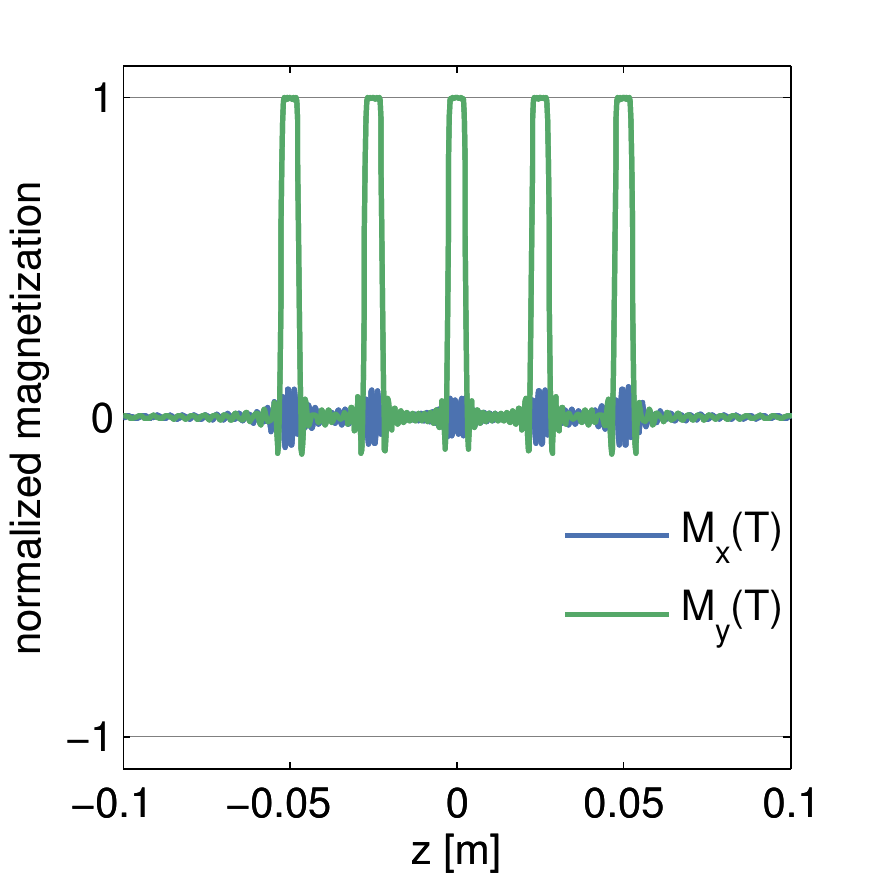}
            \hfill
        }	  

        \subcaptionbox{optimized pulse (phase shift) \label{fig:method1:3}}{%
            \includegraphics[width=0.44\textwidth]{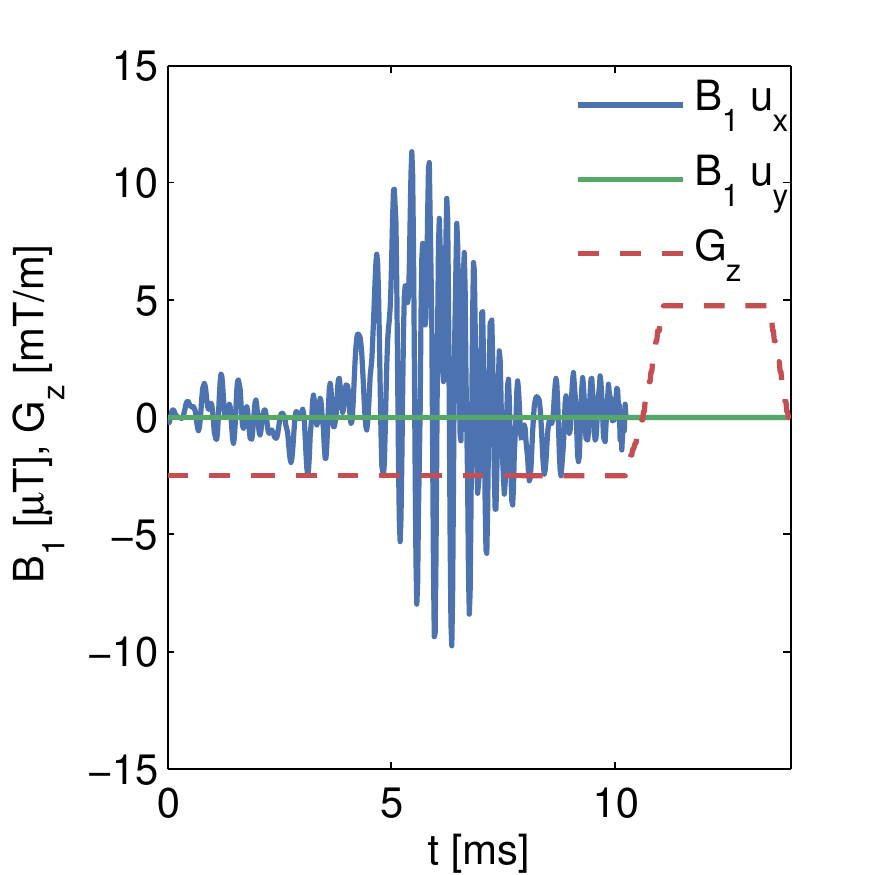}
        }
        \hfill
        \subcaptionbox{optimized slice profile of (c) \label{fig:method1:4}}{%
            \includegraphics[width=0.44\textwidth]{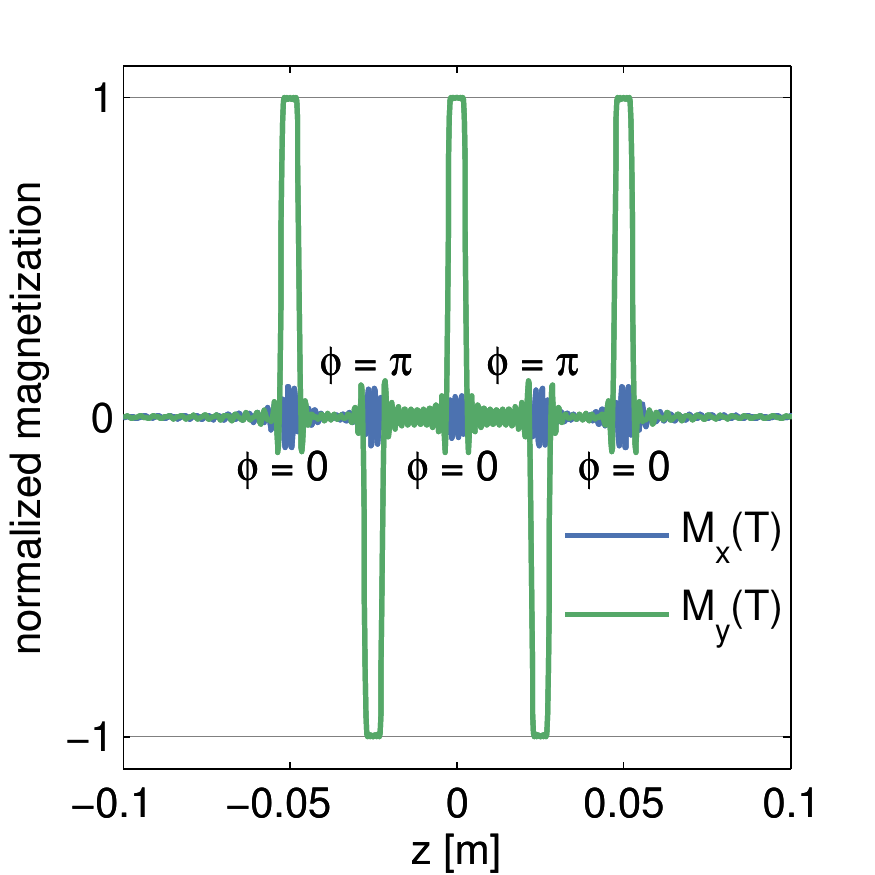}
        }
        \hfill
    \end{minipage}
    \hfill
    \begin{minipage}[c]{0.25\textwidth}
        \raggedright
        \subcaptionbox{slice-aliased reconstruction \label{fig:method1:5}}{%
            \includegraphics[width=0.98\textwidth]{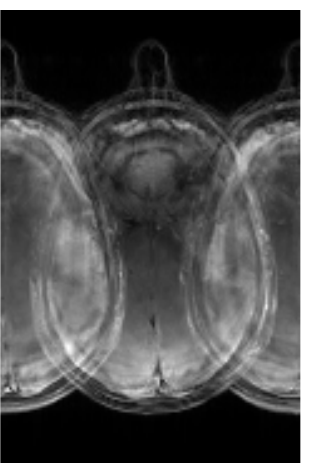} 
        }
    \end{minipage}
    \caption{Optimized pulses, slice profiles and slice-aliased Cartesian reconstruction for CAIPIRINHA-based SMS excitation pattern (five slices)} \label{fig:method1}
\end{figure}

The CAIPIRINHA-based modifications to the SMS pulse design (see \cref{fig:initial:5_zoom}) are illustrated in \cref{fig:method1} (showing the case of five slices for the sake of variation). \Cref{fig:method1:1} shows the unmodified pulse, which differs in structure from the cases with an even number of slices in, e.g., \cref{fig:sms:a6} due to the different symmetry of the slice profile (see \cref{fig:method1:2}). On the other hand, the pulse is very similar to the modified pulse for the alternating phase shift; see \cref{fig:method1:3} for the computed pulse and \cref{fig:method1:4} for the resulting slice profile. For illustration, a slice-aliased reconstruction of the acquired in-vivo data using this pulse sequence is shown in \cref{fig:method1:5}.
\begin{figure}[t]
    \centering
    \begin{minipage}[c]{0.12\textwidth}
        \centering
        z readout

        \includegraphics[width=1\textwidth,clip,trim=0 3 0 0]{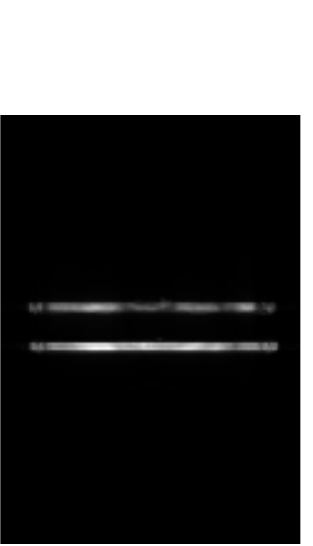}\\
        \includegraphics[width=1\textwidth,clip,trim=0 3 0 0]{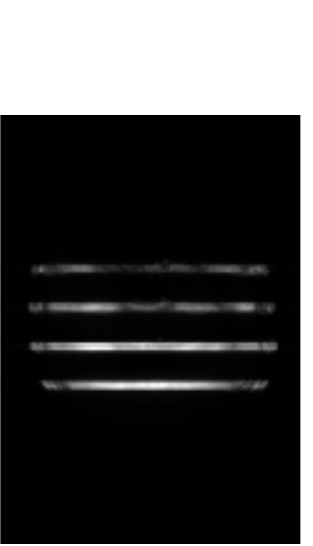}\\
        \includegraphics[width=1\textwidth,clip,trim=0 3 0 0]{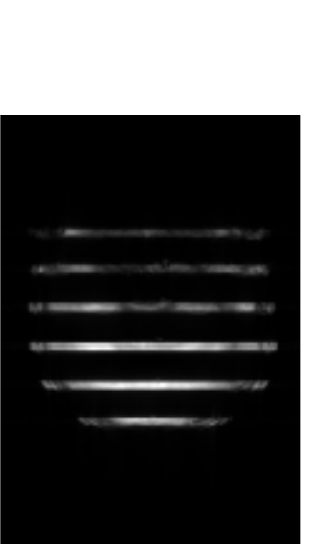}
    \end{minipage}
    \hfill
    \begin{minipage}[c]{0.8\textwidth}
        \centering
        sG reconstruction

        \includegraphics[width=0.15\textwidth]{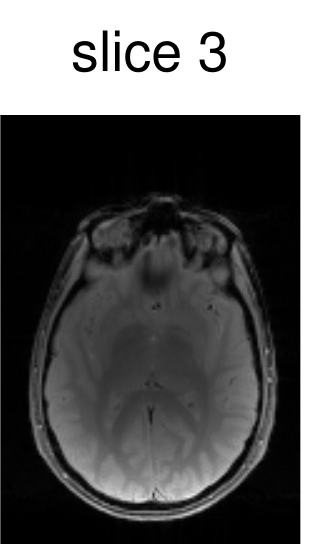}~
        \includegraphics[width=0.15\textwidth]{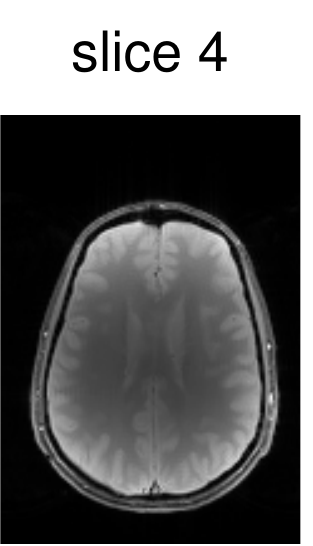}\\
        \includegraphics[width=0.15\textwidth]{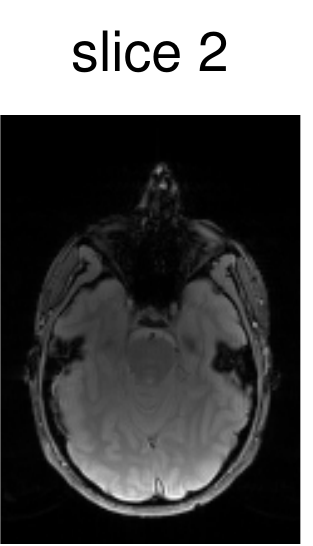}~
        \includegraphics[width=0.15\textwidth]{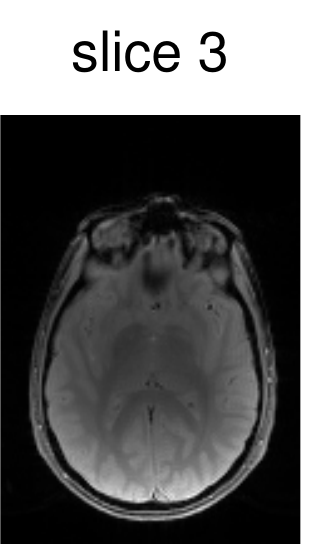}~ 
        \includegraphics[width=0.15\textwidth]{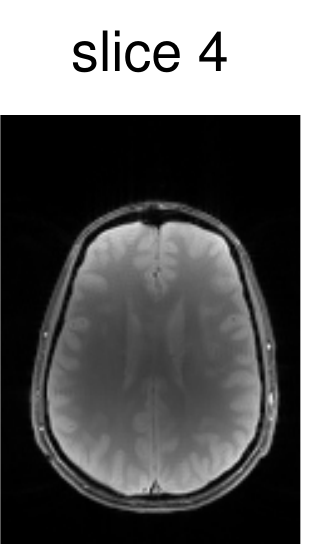}~ 	
        \includegraphics[width=0.15\textwidth]{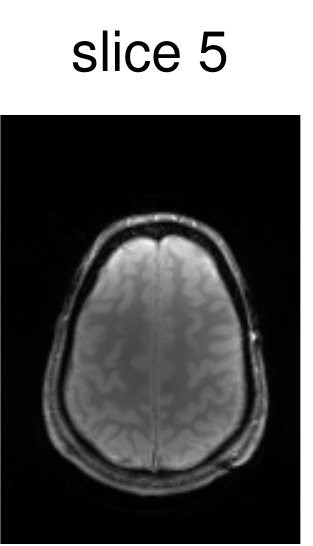}\\
        \includegraphics[width=0.15\textwidth]{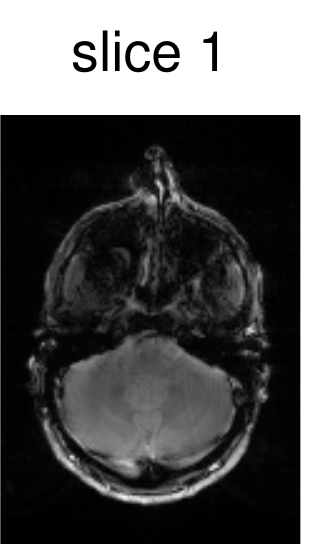}~ 
        \includegraphics[width=0.15\textwidth]{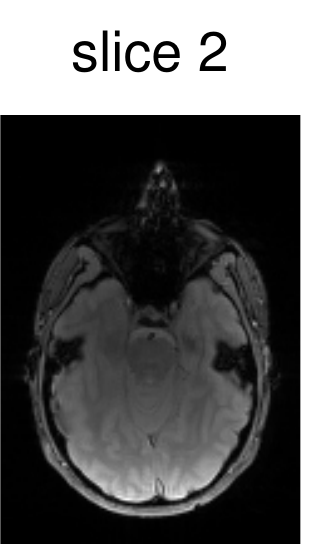}~ 
        \includegraphics[width=0.15\textwidth]{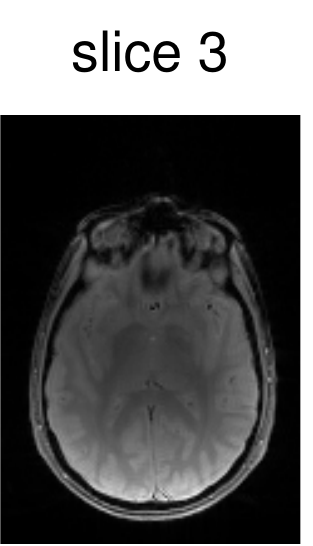}~ 
        \includegraphics[width=0.15\textwidth]{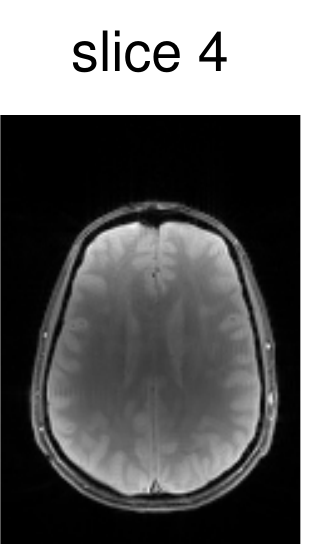}~ 
        \includegraphics[width=0.15\textwidth]{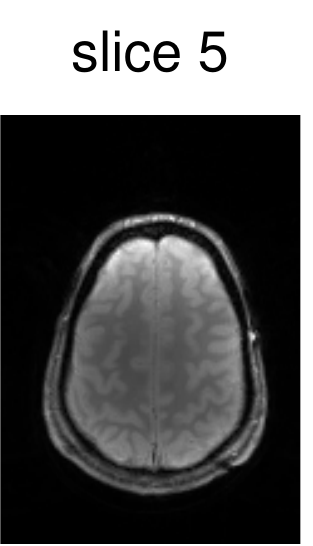}~
        \includegraphics[width=0.15\textwidth]{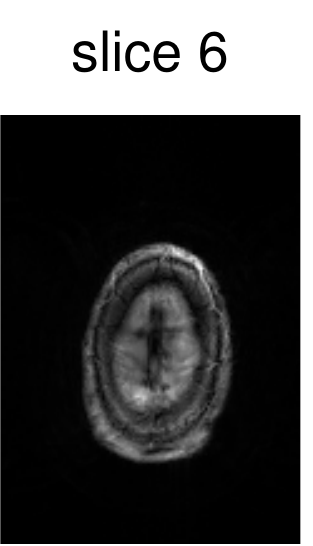} 
    \end{minipage}
    \caption{Slice-GRAPPA reconstruction of in-vivo data using CAIPIRINHA-based SMS excitation pattern for two (top), four (middle) and six (bottom) slices (left: conventional reconstruction showing the collapsed data in slice-encoding direction; right: reconstruction of GRAPPA-separated slices)}\label{fig:exp3}
\end{figure}

\Cref{fig:exp3} shows the image reconstruction using optimized RF pulses for simultaneous excitation of two, four and six slices with the same slice separation and thickness as above. As can be seen clearly in the first column, all three pulses lead to the desired excitation pattern in-vivo as well. 
The remaining columns show the slice-GRAPPA reconstructions, which illustrate that the excitation is uniform across the field of view.

\section{Discussion}
\label{sec:discussion}

Our optimization approach is related to the basic ideas presented by Conolly et al. \cite{Conolly86}. 
In the context of MRI, the implementation of this principle was also carried out by other groups using gradient \cite{Xu08, Grissom09} and quasi-Newton \cite{Vinding12} methods.
However, these methods do not make full use of second-order information and therefore achieve at best superlinear convergence.
In contrast, our Newton method makes use of exact second derivatives and is therefore quadratically convergent. 
In particular, the main contribution of our work is the efficient computation of exact Hessian actions using the adjoint approach and its implementation in a matrix-free trust-region CG--Newton method.
The use of exact derivatives speeds up convergence of the CG method, while the trust-region framework guarantees global convergence and terminates the CG method early especially at the beginning of the optimization.
Both techniques save CG steps and therefore computations of Hessian actions, allowing the use of second-order information with limited computational effort and memory requirements. 
Since computing a Hessian action incurs the same computational cost as a gradient evaluation (i.e., the solution of two ODEs; compare \eqref{eq:gradient} with \eqref{eq:Hessian}), we were able to compute a minimizer, e.g., for the single-slice example, with a computational effort corresponding to $32$ gradient evaluations ($4$ for the right-hand side in each Newton iteration and $28$ for the Hessian action in each CG iteration). This is less than the same number of iterations of a gradient or quasi-Newton method with line search (required in this case for global convergence), demonstrating the efficiency of the proposed approach.
Therefore, our method can be used to compute RF pulses with a high temporal resolution, allowing the design of pulses for a desired magnetization on a very fine spatial scale, in particular for the excitation of a sharp slice profile.

Furthermore, 
the proposed algorithm does not require an educated initial guess for global convergence (to a local minimizer, which might depend on the initial guess if more than one exists)
and allows for pulse optimization in non-standard situations where no analytic RF pulse exists (e.g., for large flip angles).
Compared to design methods using a simplification or approximation of the Bloch equation \cite{Pauly89, Pauly91}, our OC based approach is capable of including relaxation terms. However, for standard in-vivo imaging applications of the human head, the relevant relaxation times are very long compared to the RF pulse duration. 
Thus, in our examples the influence of relaxation during excitation on the designed pulses is insignificant and has been neglected in the optimization process (although the inclusion may be indicated for other applications).
The presented direct design approach allows to specify the desired magnetization in x-, y- and z-direction independently for every point in the field of view. This spatial independence of each control point allows to directly apply parallel computing to speed up the optimization process. 
While real-time optimization was not the aim of this work, a proof-of-concept implementation of the proposed approach on a GPU system (CUDA, double precision,  GeForce GTX \num{550} Ti with \num{192} cores and \SI{1024}{\mega\byte} of RAM) shows an average speedup of \num{135} (e.g., \SI{6.8}{\second} instead of \SI{989}{\second} for the single-slice example) while yielding identical results, thus making patient-specific design feasible as well as making the gap between OC and SLR pulse design nearly negligible.
This allows efficient and fast generation of accurate slice profiles -- important for minimal slice gaps, optimal contrast and low systematic errors in quantitative imaging -- for arbitrary flip angles and even for specialized pulses such as refocusing or inversion.

In particular, our approach can be used to design pulses for the simultaneous excitation of multiple slices, which increases the temporal efficiency of advanced imaging techniques such as diffusion tensor imaging, functional imaging or dynamic scans. In these contexts, SMS excitation is successfully used to reduce the total imaging time \cite{Larkman01,Nunes06,Feinberg13, Ugurbil13}; however, the peak B$_1$ amplitude of conventional SMS pulses is one of the main restrictions of applying SMS imaging to high-field systems \cite{Ugurbil13}. 
The performed studies show that compared to conventional SMS design, the presented procedure yields pulses with a reduced B$_1$ peak amplitude (e.g., \SI{40}{\percent} reduction for six simultaneous slices).
Depending on the desired temporal resolution, the bandwidth and the slice profiles of the outer slices are slightly changed, which results in a decreased B$_1$ peak amplitude. It could be shown that the peak B$_1$ amplitude does not increase linearly with the number of slices, while the power requirement per slice remains constant and the overall power consumption is comparable to that of conventional pulses. To further reduce the SAR it is necessary to either change the excitation velocity using a time-varying slice selective gradient \cite{Conolly88}, or to extend the pulse design to parallel transmit \cite{Katscher08, Wu13, Poser14, Guerin14}.  
Furthermore, our OC-based pulses produce sharp slice profiles with a lower mean absolute error compared to 
the used PM-based SLR
pulse, both in- and out-of-slice, at the cost of slightly larger out-of-slice ripples close to the in-slice regions. 
Of course, the ripple behavior of the SLR pulse can be balanced with the transition steepness by using different digital filter design methods (i.e. PM for minimizing the maximum ripple or a least squares linear-phase FIR filter for minimizing integrated squared error).
The OC ripple amplitude close to the transition band can be further controlled by using offset-dependent weights as demonstrated by Skinner et al. \cite{Skinner12_2}.
In addition, the computational complexity of OC methods is significantly higher than for direct or linearized methods. This implies that OC-based pulse design is advantageous in situations where high in-slice contrast and low B$_1$ peak amplitude are important, while SLR pulses should be used when minimal near-slice excitation and computational effort are crucial.

The presented OC approach is able to avoid some possible disadvantages of previously proposed design methods for SMS excitation. In particular, the OC design method prescribes each slice with the same uniform echo-time and phase in comparison to time-shifted \cite{Auerbach13}, phase relaxation \cite{Wong12} and nonlinear phase design techniques \cite{Zhu14, Sharma15}.
On the other hand, some of their features such as different echo times \cite{Auerbach13} or a non-uniform phase pattern \cite{Wong12, Zhu14, Sharma15} (e.g., for spin echo experiments) can be incorporated in our approach to further reduce the B$_1$ peak amplitude. 
It also should be possible to combine the OC design method with other techniques analogous to MultiPINS \cite{Eichner14, Gagoski15} that combine PINS with conventional multiband pulses for a further reduction of SAR. 
Finally, the phantom and in-vivo experiments demonstrate that it is possible to simply replace standard pulses by optimized pulses in existing imaging sequences, and that the proposed method is therefore well suited for application in a wide range of imaging situations in~MRI.

\section{Conclusions} 
\label{sec:conclusions}

This paper demonstrates a novel general-purpose implementation of RF pulse optimization based on the full time-dependent Bloch equation and a highly efficient second-order optimization technique assuring global convergence  
to a local minimizer,
which allows large-scale optimization with flexible problem-specific constraints. 
The power and applicability of this technique was demonstrated for SMS, where a reduced B$_1$ peak amplitude allows exciting a higher number of simultaneous slices or achieving a higher flip angle.
Phantom and in-vivo measurements (on a \SI{3}{\tesla} scanner) verified these findings for optimized single- and multi-slice pulses. Even for a large number of simultaneously acquired slices, the reconstructed images show good image quality and thus the applicability of the optimized RF pulses for practical imaging applications. 
While the computational requirements for optimal control approaches are of course significantly greater than for, e.g., SLR-based approaches, a proof-of-concept GPU implementation indicates that this gap can be sufficiently narrowed to make patient-specific design feasible.

Due to the flexibility of the optimal control formulation and the efficiency of our optimization strategy, it is possible to consider field inhomogeneities (B$_1$, B$_0$), design complex RF pulses for parallel transmit, or to extend the framework to include pointwise constraints due to hardware limits such as peak B1 amplitude and slew rate.

\section*{Acknowledgments}

This work is funded and supported by the Austrian Science Fund (FWF) in the context of project "SFB F3209-18" (Mathematical Optimization and Applications in Biomedical Sciences). Support from BioTechMed Graz and NAWI Graz is gratefully acknowledged. 
We would like to thank Markus Bödenler from the Graz University of Technology for the CUDA implementation of our algorithm.  

\appendix

\section{Trust-region algorithm}\label{sec:algo}

\begin{algorithm2e}[H]
    \DontPrintSemicolon
    \caption{Trust-region CG-Newton algorithm}\label{alg:trcgn}
    \KwIn{Trust region parameters $\mathrm{tol}_N$, $\mathrm{maxit}_N$, $\mathrm{tol}_C$, $\mathrm{maxit}_C$, $\rho_0$, $\rho_{\max}$, $q>1$, $0<\sigma_1<\sigma_2<\sigma_3<1$}
    \KwOut{Control ${u}$}

    Set ${u}^0\equiv 0, \quad k=0, \quad {g}\equiv 1,\quad \rho = \rho_0$\tcp*{initialization}
    \While(\tcp*[f]{TR-Newton step}){$\|g\|>\mathrm{tol}_N$ and $k<\mathrm{maxit}_N$}{
        Compute gradient ${g}({u}^k)$\;
        Set ${p}^0={r}^0=-{g}({u}^k)$, $\delta {u}=0$, $i = 0$\;
        \While(\tcp*[f]{TR-CG step}){$\|{r}^i\| > \mathrm{tol}_C\|{r}^0\|$ and $i<\mathrm{maxit}_C$}{
            Compute ${H}({u}^k){p}^i$\;
            \If(\tcp*[f]{negative curvature: CG fails}){$\langle{p}^i,{H}({u}^k){p}^i\rangle<\varepsilon$}{
                Compute $\max\{\tau:\|\delta {u} + \tau {p}^i\|\leq \rho\}$\tcp*{go to boundary of trust region}
                Set $\delta {u} = \delta {u} + \tau {p}^i$;
                \textbf{break}\;

            }
            Compute $\alpha = \|{r}^i\|/\langle {p}^i, {H}({u}^k){p}^i\rangle$\;
            \If(\tcp*[f]{step too large: model not trusted}){$\|\delta {u}+\alpha {p}^i\|\geq \rho$}{ 
                Compute $\max\{\tau:\|\delta {u} + \tau {p}^i\|\leq \rho\}$\tcp*{go to boundary of trust region}
                Set $\delta {u} = \delta {u} + \tau {p}^i$;
                \textbf{break}\;
            }
            Set ${r}^{i+1} = {r}^{i} - \alpha {H}({u}^k){p}^i$\;
            Set ${p}^{i+1} = {r}^{i+1} + \|{r}^{i+1}\|^2/\|{r}^i\|^2 {p}^{i}$\;
            Set $\delta {u} = \delta {u} + \alpha {p}^i$,\quad $i = i+1$\; 
        }
        Compute $\delta J_a = J({u}^k) - J({u}^k +\delta {u})$\tcp*{actual function decrease}
        Compute $\delta J_m = -\tfrac12 \langle \delta {u},{H}({u}^k)\delta {u}\rangle - \langle\delta {u}, {g}({u}^k)\rangle$\tcp*{predicted function decrease}
        \If(\tcp*[f]{sufficient decrease}){$\delta J_a >\varepsilon$ and $\delta J_a > \sigma_1 \delta J_m$}{
            Set ${u}^{k+1} = {u}^k +\delta {u}$\tcp*{accept step}
        }
        \uIf(\tcp*[f]{step accepted, model good}){$\delta J_a>\varepsilon$ and $|\delta J_a /\delta J_m - 1| \leq 1-\sigma_3$}{
            Set $\rho = \min\left\{q \rho,\rho_{\max}\right\}$\tcp*{increase radius}
        }
        \ElseIf(\tcp*[f]{step rejected, no decrease}){$\delta J_a \leq \varepsilon$}{
            Set $\rho = \rho / q$\tcp*{decrease radius}
        }
        \ElseIf(\tcp*[f]{model bad}){$\delta J_a < \sigma_2 \delta J_m$}{
            Set $\rho = \rho / q$\tcp*{decrease radius}
        }
    }
\end{algorithm2e}

\newpage
\section{Discretization}\label{sec:discrete}

\paragraph{Cost functional:}
\[
    J({M},{u})= \frac{1}{2} \sum_{i=1}^Z \Delta z_i |{M}_{N,i}-{M_d}(z_i)|_2^2 + \frac\alpha2 \sum_{m=1}^N \Delta t_m |{u}_m|_2^2
\]
\paragraph{Bloch equation} for all $i=1,\dots,Z$:
\[
    \begin{aligned}[t]
        \left[{I}-\frac{\Delta t_m}{2}{A}({u}_m;z_i)\right] {M}_{m,i} &= \left[{I}+\frac{\Delta t_m}{2}{A}({u}_m;z_i)\right] {M}_{m-1,i} +\Delta t_m {b}, \quad m=1,\dots,N\\
        {M}_{0,i} &= {M^0}(z_i)
    \end{aligned}
\]
\paragraph{Adjoint equation} for all $i=1,\dots,Z$:
\[
    \begin{aligned}[t]
        \left[{I}-\frac{\Delta t_m}{2}{A}({u}_m;z_i)^T\right] {P}_{m,i} &= \left[{I}+\frac{\Delta t_{m+1}}{2}{A}({u}_{m+1};z_i)^T\right] {P}_{m+1,i} , \quad m=1,\dots,N-1\\
        \left[{I}-\frac{\Delta t_N}{2}{A}({u}_N;z_i)^T\right]{P}_{N,i}&={M}_{N,i}-{M_d}(z_i)\\
    \end{aligned}
\]
\paragraph{Discrete gradient} for all $m=1,\dots,N_u$: \quad $\bar{{M}}_{m} := \frac12({M}_m + {M}_{m-1})$,
\[
    {g}_m = \alpha {u}_m + \gamma B_1
    \begin{pmatrix}
        \sum_{i=1}^Z \Delta z_i \left({P}_{m,i}^T {A_1} \bar{{M}}_{m,i}\right) \\ 
        \sum_{i=1}^Z \Delta z_i \left({P}_{m,i}^T {A_2} \bar{{M}}_{m,i}\right)
    \end{pmatrix}
\]
\paragraph{Linearized state equation} for all $i=1,\dots,Z$:
\[
    \begin{aligned}[t]
        \left[{I}-\frac{\Delta t_m}{2}{A}(u_m;z_i)\right] \delta {M}_{m,i} &= \left[{I}+\frac{\Delta t_m}{2}{A}(u_m;z_i)\right] \delta {M}_{m-1,i} + \Delta t_m {A'}(\delta u_m) \bar{{M}}_{m,i} , \quad m=1,\dots,N\\
        \delta {M}_{0,i} &= 0
    \end{aligned}
\]
\paragraph{Linearized adjoint equation} for all $i=1,\dots,Z$:
\[
    \begin{aligned}[t]
        \left[{I}-\frac{\Delta t_m}{2}{A}(u_m;z_i)^T\right] \delta {P}_{m,i} &= \left[{I}+\frac{\Delta t_{m+1}}{2}{A}(u_{m+1};z_i)^T\right] \delta {P}_{m+1,i} +\frac{\Delta t_m}{2}{A'}(\delta {u}_m)^TP_{m,i}\\
        \MoveEqLeft[-1]+\frac{\Delta t_{m+1}}{2}{A'}(\delta {u}_{m+1})^T {P}_{m+1,i} , \quad m=1,\dots,N-1\\
        \left[{I}-\frac{\Delta t_N}{2}{A}(u_N;z_i)^T\right]\delta {P}_{N,i}&=\delta {M}_{N,i} + \frac{\Delta t_N}{2}{A'}(\delta {u}_{N})^T{P}_{N,i}
    \end{aligned}
\]
\paragraph{Discrete Hessian action} for all $m=1,\dots,N_u$: \quad $\bar {\delta {M}_{m}} := \frac12(\delta {M}_m + \delta {M}_{m-1})$
\[
    \quad [{H}({u}){h}]_m= \alpha {h}_{m} + \gamma B_1
    \begin{pmatrix}
        \sum_{i=1}^Z \Delta z_i \left(\delta {P}_{m,i}^T {A_1} {\bar M}_{m,i} + {P}_{m,i}^T {A_1} {\bar{\delta M}}_{m,i}\right)\\ 
        \sum_{i=1}^Z \Delta z_i \left(\delta {P}_{m,i}^T {A_2} {\bar M}_{m,i} + {P}_{m,i}^T {A_2} {\bar{\delta M}}_{m,i}\right)
    \end{pmatrix} 
\]

\clearpage

\printbibliography
\end{document}